\documentclass[12pt,reqno]{amsart}
\usepackage{amsthm}
\usepackage{fullpage,url}
\usepackage[dvips]{hyperref}
\usepackage{amscd}   
\usepackage[all,graph,poly]{xy} 
\usepackage{amssymb,amsmath}
\usepackage{graphicx}
\usepackage{pstricks,pst-plot,pst-node}
\usepackage[dvips]{hyperref}
\usepackage{multicol}
\usepackage{setspace}
\DeclareFontEncoding{OT2}{}{} 

\newcommand{\abs}[1]{\lvert#1\rvert}

\newcommand{\braid}{\between}
\newcommand{\CC}{\mathbb{C}}
\newcommand{\cE}{\mathcal{E}}
\newcommand{\cG}{\mathcal{G}}

\newcommand{\cM}{\mathcal{M}}
\newcommand{\cO}{\mathcal{O}}
\newcommand{\comm}{\perp }

\newcommand{\ip}[2]{\langle #1, #2 \rangle}
\newcommand{\op}[1]{\operatorname{#1}}
\newcommand{\PP}{\mathbb{P}}
\newcommand{\QQ}{\mathbb{Q}}
\newcommand{\rel}[2]{\mathcal{C}_{#1} \langle #2 \rangle }

\newcommand{\RR}{\mathbb{R}}

\newcommand{\ZZ}{\mathbb{Z}}
\DeclareMathOperator{\Aut}{Aut}

\DeclareMathOperator{\sym}{S}
\newcommand{\e}[1]{\ar@{-}[#1]}
\newcommand{\ed}[1]{\ar@{--}[#1]}
\newcommand{\dE}[1]{\ar@{=}[#1]}
\newcommand{\er}[1]{\ar@{-}[#1] |-{\SelectTips{cm}{}\object@{<}} |{\SelectTips{eu}{}\object@{}}}
\newcommand{\el}[1]{\ar@{-}[#1] |-{\SelectTips{cm}{}\object@{>}} |{\SelectTips{eu}{}\object@{}}}
\newcommand{\edr}[1]{\ar@{--}[#1] |-{\SelectTips{cm}{}\object@{<}} |{\SelectTips{eu}{}\object@{}}}
\newcommand{\edl}[1]{\ar@{--}[#1] |-{\SelectTips{cm}{}\object@{>}} |{\SelectTips{eu}{}\object@{}}}
\newcommand{\edb}[1]{\ar@{=}[#1] |-{\SelectTips{cm}{}\object@{<}} |{\SelectTips{eu}{}\object@{}}}
\newcommand{\edbd}[1]{\ar@{==}[#1] |-{\SelectTips{cm}{}\object@{<}} |{\SelectTips{eu}{}\object@{}}}

 \newcommand{\lul}[1]{\ar@{}[l]_<<{#1}}
\newcommand{\rrul}[1]{\ar@{}[r]^<<<<{#1}}
\newcommand{\rul}[1]{\ar@{}[r]^<<{#1}}
\newcommand{\ldl}[1]{\ar@{}[l]^<<{#1}}
\newcommand{\rdl}[1]{\ar@{}[r]_<<{#1}}
\newcommand{\dl}[1]{\ar@{}[d]_<<{#1}}
\newcommand{\dll}[1]{\ar@{}[dd]_{#1}}

\providecommand{\abs}[1]{\lvert#1\rvert}
\swapnumbers
\newtheorem{theorem}{Theorem}[section]
\newtheorem{lemma}[theorem]{Lemma}

\newtheorem{proposition}[theorem]{Proposition}

\theoremstyle{definition}
\newtheorem{definition}[theorem]{Definition}

\newtheorem{topic}[theorem]{}

\theoremstyle{remark}
\newtheorem{remark}[theorem]{Remark}
\newtheoremstyle{head}
{}
{}
{\bfseries}
{}
{}
{}
{.5em}
{}

\theoremstyle{head}

\begin{document}
%
%
{\bf {\large On Coxeter diagrams of complex reflection groups}}
\par
\vspace{.5cm}
Author : {\large Tathagata Basak}\par
\vspace{.5cm}
{\small Address : Iowa State University, Department of mathematics, Carver Hall, Ames, IA 50011.}
\par
{\small email : tathastu@gmail.com}
\newline
\newline
{\small 
Abstract: We study Coxeter diagrams of some unitary reflection groups.
Using solely the combinatorics of diagrams, we give a new proof of the
classification of root lattices defined over $\cE = \ZZ[e^{ 2 \pi i/3}]$:
there are only four such lattices, namely, the $\cE$--lattices whose real
forms are $A_2$, $D_4$, $E_6$ and $E_8$.
Next, we address the issue of characterizing the diagrams for unitary
reflection groups, a question that was raised by Brou\'{e}, Malle and Rouquier.
To this end, we describe an algorithm which,  given a unitary reflection
group $G$, picks out a set of complex reflections. The algorithm
is based on an analogy with Weyl groups. If $G$ is a Weyl
group, the algorithm immediately yields a set of simple roots. 
Experimentally we observe that if $G$ is primitive and $G$ has a set of roots
whose $\ZZ$--span is a discrete subset of the ambient vector space,
then the algorithm selects a minimal generating set for $G$.
The group $G$ has a presentation on these generators such that
if we forget that the generators have finite order then we get
a (Coxeter-like) presentation of the corresponding braid group.
For some groups, such as
$G_{33}$ and $G_{34}$, new diagrams are obtained.
For $G_{34}$, our new diagram 
extends to an ``affine diagram" with $\ZZ/7\ZZ$ symmetry. 
\newline
\newline
{\it Keywords: Unitary reflection group, Coxeter diagram, Weyl group, simple root.} \\
2000 {\it Mathematics subject classification:} 20F55, 20F05, 20F65, 51F25.
}
%
%
\section{Introduction}
%
%
\begin{topic}{\bf Background on unitary reflection groups:}
\label{background}
 Let  $G$ be a finite subgroup of the unitary group $U(n)$ generated by complex reflections, such that
$G$ acts irreducibly on $\CC^n$. We shall simply say that $G$ is an unitary reflection group.
Let $\cM$ be the union of the fixed point sets of the complex reflections in $G$
and let $X_G = \CC^n \setminus \cM$. The fundamental group $\op{Braid}(G) := \pi_1( X_G/G)$ is called
the {\it (generalized) braid group associated to $G$}.
Let $k$ be the minimum number of complex reflections needed to generate $G$.
We say that $G$ is {\it well generated} if $k = n$.
The smallest subfield $F$ of $\CC$ that contains all the complex character
values of $G$, is called the {\it field of definition} of $G$. We say that $G$ is
defined over $F$.
\par
Unitary reflection groups were classified by
Shephard and Todd in \cite{st:fu}. For a self-contained proof of the
classification, which is similar in spirit to part of our work,
see \cite{amc:fc}. 
A convenient table of all these groups and their properties may
be found in \cite{bmr:cr}. There is an infinite series, denoted by
$G(de, e,r)$, and $34$ others, denoted by $G_4, G_5, \dotsb, G_{37}$.
Unitary reflection groups have many invariant theoretic properties that
are similar to those of the orthogonal
reflection groups. Most of these properties were initially established
for the unitary reflection groups, via case by case verification through 
Shephard and Todd's list. Recently, there has been a lot of progress in trying
to find unified and more conceptual proofs. (For example, see
\cite{db:fc} and the references therein). However, a coherent theory,
like that of the classical Coxeter groups and Weyl groups, is still
not in place and many mysteries still remain. One of these mysteries
involve the diagrams for unitary reflection groups. 
\par
Coxeter presentations of orthogonal reflection groups are encoded in their
Coxeter-Dynkin diagrams. Similarly, for each unitary reflection group $G$, there is a diagram
$D_G$, that encodes a presentation of $G$ (Such a $D_G$ is given in \cite{bmr:cr}
for all but six groups. For the remaining six groups presentations of the
corresponding braid group were conjectured in \cite{bm:ep}, in terms of certain
diagrams. The proof of this conjecture was completed in \cite{db:fc}).
Most of these diagrams were first introduced by Coxeter in \cite{hc:fu}. The vertices of
$D_G$ correspond to complex reflections that form a minimal set of generators for $G$.
Other than that, the definition of $D_G$ is ad-hoc and case by case.
It is curious that even though these diagrams do not have any uniform
definition, they contain a lot of non-trivial information
about the group $G$. We quote two sample results which state that the 
weak homotopy type of $X_G/G$ and the invariant degrees of $G$ can be recovered from the
diagram $D_G$.
\begin{theorem}[\cite{db:fc}]  (a) The universal cover of $X_G/G$ is contractible.
\par
(b) $G$ has a minimal generating set of complex reflections, $\bar{\mathbf{b}} = \lbrace \bar{b}_1, \dotsb, \bar{b}_k \rbrace$,
which can be lifted to a set of generators $\mathbf{b} = \lbrace b_1, \dotsb, b_k \rbrace$ of $\op{Braid}(G)$, with the following
property: There is a set of positive homogeneous relations $R'_G(\mathbf{b})$ in the alphabet $\mathbf{b}$
such that $\op{Braid}(G)$ and $G$ have the following presentations:
\begin{equation*} 
\op{Braid}(G) \simeq \langle \mathbf{b} \vert R'_G(\mathbf{b}) \rangle, \text{\; \;}
G \simeq \langle \bar{\mathbf{b}} \vert R'_G(\bar{\mathbf{b}} ), \bar{b}_i^{n_i} = 1 \text{\; for all \;} i\rangle,
 \end{equation*}
where $n_i$ is the order of $\bar{b}_i$ in $G$.
(These presentations are encoded by a diagram $D_G$ with vertex set $\mathbf{b}$, 
the edges indicating the relations in $R'_G$). 
\label{t-brd}
\end{theorem}
\begin{theorem}[Th. 5.5 of  \cite{os:ur}, \cite{bm:ep}, \cite{db:fc}]  Let $G$ be a {\it well generated}
unitary reflection group. The {\it Coxeter number} of $G$ is defined to be the largest positive
integer $h$ such that $e^{2 \pi i/h}$ is an eigenvalue of an element of $G$.
Then the product of the generators of $G$ corresponding to the vertices of $D_G$, in
certain order, has eigenvalues $e^{ 2 \pi i(d_j-1)/h}$ where $d_j$ are the
invariant degrees of $G$.
\label{t-cox}
\end{theorem}
If $G$ is a Weyl group and $V$ is the complexification of the standard representation of $G$, then
both \ref{t-brd} (due to Briskorn, Saito, Deligne) and \ref{t-cox} (probably due to
Borel, Chevalley, Steinberg) are classical.
These results were first verified for most of the groups in Shephard and
Todd's list by arguments split into many separate cases. Essentially classification
free proofs are now known, by recent work of David Bessis (see \cite{db:fc}, 
where the long standing task of showing that $X_G$ is $K(\pi, 1)$ and finding
Coxeter-like presentation for $\op{Braid}(G)$ were completed).
But we still do not know of a way to characterize the diagrams for unitary
reflection groups.
\par
In this article, we study the diagrams of a few unitary reflection groups.
The main results are discussed below. They are motivated by analogies with Weyl groups.
\end{topic}
\begin{topic}{\bf Summary of results: }
Our approach is to view unitary reflection groups as sets of automorphisms
of ``complex lattices''.
Let $\cE = \ZZ[e^{ 2 \pi i/3}]$. The main examples of unitary reflection groups,
that we want to study, act as automorphisms of a sequence of $\cE$--lattices, namely,
$ A_2^{\cE} \subseteq D_4^{\cE} \subseteq E_6^{\cE} \subseteq E_8^{\cE}$.
Our interest in these lattices stems from their importance in studying the
complex hyperbolic reflection group with $Y_{555}$ diagram and its conjectured
connection with the bimonster (see \cite{dja:Y_5552}, \cite{tb:el}).
In section \ref{seceirootlattice}, we present a new proof of theorem 2.2 of 
\cite{dja:Y_5552}, which states that $ A_2^{\cE}, D_4^{\cE}, E_6^{\cE}, E_8^{\cE}$
are the only  ``$\cE$--root lattices".  Our proof is like the $A$--$D$--$E$
classification of Euclidean root lattices and is similar in spirit to the arguments
in \cite{amc:fc} and \cite{mh:cr}. It is purely a linear algebra argument
that only uses the diagrams for the complex reflection groups.
This proof should be viewed as an illustration of the usefulness of
the ``complex diagrams". 
\par
In Section \ref{secweylvector} 
we address the following question, that was raised in \cite{bmr:cr}: 
{\it How to characterize the diagrams for the unitary reflection groups?}
To this end, we describe an algorithm (in \ref{method} and \ref{algo})
which, given the group $G$, the integer $k$ and a random vector in $V$, selects a set  
$\bar{\mathbf{a}} = \lbrace \bar{a}_1, \dotsb, \bar{a}_k \rbrace$
of reflections in $G$. Our algorithm is based on a
generalization of a ``Weyl vector''. We show that ``Weyl vectors'' exist
for all unitary reflection groups (see theorem \ref{prop-fixpt}).
If $G$ is a Weyl group, then one can easily check that
$\bar{\mathbf{a}}$ is a set of simple roots of $G$.
\par
If $G$ is primitive and defined over an imaginary quadratic
extension of $\QQ$, then we experimentally observe that $\bar{\mathbf{a}}$
is a minimal set of generators of $G$.
There exists a set of positive homogeneous relations 
$R_G(\mathbf{a})$ in the alphabet $\mathbf{a} = \lbrace a_1, \dotsb, a_k \rbrace$
such that:
\newline
\newline
{\it In every execution of the algorithm, the generators $\bar{\mathbf{a}} $
satisfy the relations $R_G(\bar{\mathbf{a}})$.}
\newline
\par
We find that the reflections $\bar{\mathbf{a}}$ form Coxeter's diagram in the examples
of our main interest, namely, the reflection groups related to the $\cE$--root lattices.
For some groups $G$, namely, $G_{12}$, $G_{29}$, $G_{31}$, $G_{33}$ and $G_{34}$,
new diagrams are obtained. More precisely, the generators
$\bar{\mathbf{a}}$ selected by the algorithm \ref{algo} do not satisfy the relations
$R'_G$ known from \cite{bmr:cr}, \cite{bm:ep}.
In section \ref{s-newd}, we verify that: 
\newline
\newline
{\it   The group $\op{Braid}(G)$ has a presentation given by
$\langle \mathbf{a} \vert R_G( \mathbf{a}) \rangle$ and $G$ has a presentation given by 
$\langle \bar{\mathbf{a}} \vert R_G( \bar{\mathbf{a}}), \bar{a}_i^{n_i} = 1 \text{\; for all \;} i \rangle$,
that is, Theorem \ref{t-brd}(b) holds for the new diagrams (see \ref{p-mutate}). }
\newline
\par
We have verified that Theorem \ref{t-cox} also holds for the new diagrams for 
$G_{29}$, $G_{33}$ and $G_{34}$. The other two groups $G_{12}$ and $G_{31}$ are not well-generated.
\par
Let $G \in \lbrace G_{29}, G_{31}, G_{33}, G_{34} \rbrace$.
For these groups, the relations $R_G$ that our generators satisfy are different from those
previously known.
We note that all the relations in $R_G$ (i.e. those needed to present $\op{Braid}(G)$)
are of the form $x_0 x_1 \dotsb x_{m-1} = x_1 x_2 \dotsb x_{m}$, for a set of generators
$\lbrace x_i \colon  i \in \ZZ/k \ZZ \rbrace$, which form
a minimal cycle in the diagram. When $k = 2$, these are the
Coxeter relations. For most $G$, the group $\op{Braid}(G)$
has a presentation consisting of only this kind of relations (see the table in \cite{bmr:cr}).
Following Conway \cite{ccs:26}, we call these deflation relations.
The deflation relations encountered in $G_{33}$ and $G_{34}$ are, moreover, all ``cyclic"
(see \ref{l-cyclic}, \ref{l-rel}). 
For $G_{29}$, $G_{33}$ or $G_{34}$, a presentation
of the corresponding braid group is obtained by taking one braid relation
for each edge and one deflation relation for each minimal cycle in the diagram.
This makes us wonder if the right
notion of a diagram for these groups is the $2$ dimensional polyhedral complex
obtained by attaching $2$-cells to the minimal cycles in the graphs $D_G$, so that,
finiteness of $G$ translates into the vanishing of the first homology of the polyhedral complex.
\par
Note that $G_{33}$ and $G_{34}$, are part of a few exceptional cases, in which, the
diagrams known in the literature do not have some of the desirable properties.
(For example, see question
2.28 in \cite{bmr:cr} and the remark following it). So there seems to
be some doubt whether the diagrams known in the literature for these examples,
are the ``right ones''.
\par
In section \ref{secaffine} we describe affine diagrams
for unitary reflection groups defined over $\cE$.
The affine diagrams are obtained from the unitary
diagrams by adding an extra node. They encode presentations for the
corresponding affine complex reflection groups.
For each affine diagram, we describe a ``balanced numbering" on its vertices,
like the
$\begin{smallmatrix} & & 3& & & & & \\ 2 & 4 & 6 & 5 & 4 & 3 & 2 & 1 \end{smallmatrix} $
numbering on the affine $E_8$ diagram.
The existence of a balanced numbering on a diagram implies that the corresponding
reflection group is not finite. So an affine diagram cannot occur as
a full sub-graph of a diagram for a unitary reflection group.
These facts, coupled with a combinatorial argument, complete the classification
of $\cE$--root lattices.
The affine diagrams are often more symmetric compared
to the unitary diagrams, like in the real case. For example, for $G_{34}$, (which is the
reflection group of the Coxeter-Todd lattice
$K_{1 2}^{\cE}$), we get an affine diagram with rotational $\ZZ/7\ZZ$ symmetry.
\par
The complex root systems that we use in computer experiments are described in 
the appendix \ref{ap-1}.  Appendix \ref{ap-2} contains 
details of proofs of statements made in section \ref{s-newd}.
\end{topic}
%
%
\begin{topic}{\bf Shortcomings of algorithm \ref{algo}: } If $G$ is imprimitive or not defined over 
$\QQ$ or an imaginary quadratic extension of $\QQ$, then our algorithm does not work,
in the sense that the set $\bar{\mathbf{a}}$ of reflections chosen by the algorithm
usually do not form a minimal set of generators of $G$. 
Also, one knows from \cite{bmr:cr} and \cite{bm:ep} that
the braid groups of $G_{12}$ and $G_{24}$ can be presented using cyclic deflation
relations, but we were unable to find such presentations of these braid groups on the
generators selected by our algorithm.
So while the algorithm \ref{algo} does not provide a definite characterization
of the diagrams for the unitary reflection groups,  the observations in the previous paragraphs
seem to indicate that the diagrams have a geometric origin.
\end{topic}
%
%
We finish this section by introducing some basic definitions and notations
to be used. 
\begin{topic}
{\bf Reflection groups and root systems: }
\label{def-reflection}
Let $V$ be a complex vector space with an hermitian form
(always assumed to be linear in the second variable).
If $x \in V$, then $\abs{x}^2 = \ip{x}{x}$ is called the {\it norm} of $x$.
Given a vector $x$ of non-zero norm and a root of unity $u \neq 1$, let
\begin{equation*}
\phi_x^u(y) = y - (1 - u)\ip{x}{y}\abs{x}^{-2} x.
\end{equation*}
The automorphism $\phi_x^u$ of the hermitian vector space $V$ is
called an $u$--{\it reflection} in $x$, or simply, a {\it complex reflection}.
The hyperplane $x^{\bot}$ (or its image in the projective space $\PP(V)$),
fixed by $\phi_x^u$, is called the {\it mirror} of reflection.
A {\it complex reflection group} $G$ is a discrete subgroup of
$\op{Aut}(V, \ip{\;}{\:})$, generated by complex reflections.
A mirror of $G$ is a hyperplane fixed by a reflection in $G$.
A complex reflection (resp. complex reflection group)
is called a {\it unitary reflection} (resp. {\it unitary reflection group})
if the hermitian form on $V$ is positive definite.
We shall omit the words ``complex" or ``unitary", if they are clear from context.
A unitary reflection group $G$ acting on $V$ is {\it reducible}
(resp. {\it imprimitive}) if $V$ is a direct sum $V = V_1 \oplus \dotsb \oplus V_l$
such that $ 0 \neq V_1 \neq V$ and each $V_j$ is fixed by $G$
(resp.  the collection of $ V_j $ is stabilized by $G$).
Otherwise $G$ is {\it irreducible} (resp. {\it primitive}).
Unless otherwise stated, we always assume that $G$ is irreducible. 
\par
Let $G$ be a unitary reflection group acting on $\CC^k$ with the standard
hermitian form. Let $F$ be the field of definition of $G$.
Let $\cO$ be the ring of integers in $F$.
Let $\cO^*$ be the group of units of $\cO$.
A vector $r$ in an $\cO$--module $K$ is {\it primitive}
if $r = m r'$ with $m \in \cO$ and $r' \in K$ implies that
$m \in \cO^*$.
Let $\Phi$ be a set of primitive vectors in $\cO^k$ such that:
\begin{itemize}
\item $\Phi$ is stable under the action of $G$, 
\item $\lbrace r^{\bot} \colon r \in \Phi \rbrace$ is equal to the set
of mirrors of $G$, and
\item given $r \in \Phi$, $u r \in \Phi$ if and only if $u$ is an unit
of $\cO$.
\end{itemize}
Such a set of vectors will be called a {\it (unitary) root system}
for $G$, defined over $\cO$. The group $\cO^*$ acts on
$\Phi$ by multiplication. An orbit is called a {\it projective root}.
A set of projective roots for $G$ is denoted by $\Phi_*$ or $\Phi_*(G)$.
\end{topic}
\begin{topic}{\bf Lattices and their reflection groups: }
\label{def-lattice}
Let $F$ be a number field. Let $\cO$ be the ring of integers of $F$.
Assume that $\cO$ is a unique factorization domain.
Fix an embedding of $F$ in $\CC$ and identify $F$ and $\cO$ as
subsets of $\CC$ via this embedding. Assume that $\cO$ forms a
discrete set in $\CC$. The examples that will be important to us
are the integers,  the Gaussian integers $\cG = \ZZ[i]$
and the Eisenstein integers $\cE = \ZZ[e^{2 \pi i/3}]$.
(we also briefly consider $\ZZ[\sqrt{-2}]$ and $\ZZ[(\sqrt{-7} + 1)/2]$).
To fix ideas, one may take $\cO = \cE$.
In the next section we only work over this ring.
\par
A {\it lattice} $K$, defined over $\cO$, is a free $\cO$--module of
finite rank with an $\cO$--valued hermitian form. 
Let $V = \CC \otimes_{\cO} K $ be the complex vector space
underlying $K$. The {\it dual lattice} of $K$, denoted by $K'$, is
the set of vectors $y \in V$ such that $\ip{y}{x} \in \cO$
for all $x \in K$. 
\par
A {\it root} of $K$ is a primitive vector $r \in K$ of non-zero norm such that
$\phi_r^u \in \Aut(K)$ for some root of unity $u \neq 1$. 
The {\it reflection group} of $K$, denoted by $R(K)$, is the
subgroup of $\Aut(K)$ generated by reflections in the roots of $K$.
The projective roots of $K$ are in bijection with the mirrors of $R(K)$.
If $K$ is positive definite, then the roots of $K$ form a
unitary root system, denoted by $\Phi_K$, for the unitary reflection
group $R(K)$.
\end{topic}
\begin{topic}{\bf (Root) diagrams: }
\label{def-diagram}
Consider the permutation matrices acting on $k \times k$ hermitian matrices
by conjugation. An orbit $D$ of this action is called a {\it (root) diagram} or simply a diagram.
(There is a closely related notion of Coxeter diagram defined in section \ref{s-newd}).
If $M = ( \! ( m_{i j} ) \! )$ is a representative of an orbit $D$, then we say that
$M$ is a {\it gram matrix} of $D$.
Let $\Delta = \lbrace r_1, \dotsb, r_k \rbrace$ be a subset of a hermitian vector
space $V$. Let $m_{i j} =\ip{r_i}{r_j}$. The matrix $( \! ( m_{i j} ) \! )$ is called a {\it gram matrix}
of $\Delta$ and the corresponding diagram is denoted by $D(\Delta)$.
Let $\Phi$ be a root system for a unitary reflection group $G$.
If $\Delta$ is a subset of $\Phi$ such that
$\lbrace \phi_r^{u_r} \colon r \in \Delta \rbrace$ is a 
minimal generating set for $G$ (for some units $u_r$), then
we say that $D(\Delta)$ is a {\it root diagram} for $G$. 
\par
Pictorially, a diagram $D$ is conveniently
represented by drawing a directed graph $D$ with labeling of vertices and edges,
as follows: Let $v(D) = \lbrace x_1, \dotsb, x_k \rbrace$ be the set of vertices
of $D$. We remember the entry $m_{i i}$ by labeling the vertex $x_i$ with $m_{i i}$.
We remember the entry $m_{i j}$ by drawing a directed edge from $x_j$ to $x_i$
labeled with $m_{i j}$ or equivalently, by drawing a directed edge from $x_i$ to $x_j$
labeled with $\bar{m}_{i j}$ (but not both).
\par
Let $D$ be a diagram with gram matrix $( \! (m_{i j} ) \! )$. Assume that
$m_{i j} \in \cO$ for all $i$ and $j$. 
Define $L(D)$ to be the $\cO$--lattice generated by linearly independent
vectors $\lbrace x_1, \dotsb, x_k \rbrace$ with $\ip{x_i}{x_j} = m_{ i j}$.
Conversely, let $L$ be an $\cO$--lattice having a set of roots $\Delta$ which form a
minimal spanning set for $L$ as an $\cO$--module. Then the diagram $D(\Delta)$
is called a {\it root diagram} or simply a diagram for $L$. 
Let $D$ be a diagram for $L$.  Then $L(D)$
surjects onto $L$ preserving the hermitian form. If the gram matrix of
$D$ is positive definite, then $L(D) \simeq L$.
We shall usually denote the vertices of $D$ and the corresponding vectors
of $L$ by the same symbol.
\par
Two diagrams $D$ and $D'$ are {\it equivalent} if $L(D) \simeq L(D')$.
In this case, we write $D \simeq D'$. 
Let $L = L(D)$, $v(D) = \lbrace x_1, \dotsb, x_k \rbrace$ and $u_1, \dotsb, u_k$
be units. Then there is a diagram for $L$ whose vertices correspond to the
generators $\lbrace u_1 x_1, \dotsb, u_k x_k \rbrace$.
These two diagrams are equivalent. The only difference between them is in the edge
labeling, which may differ by units.
\end{topic}
{\bf Acknowledgments: } I would like to thank Prof. Daniel Allcock, Prof. Jon Alperin,
Prof. Michel Brou\'{e}, Prof. George Glauberman and Prof. Kyoji Saito for useful discussions and my advisor Prof.
Richard Borcherds for his help and encouragement in the early stages of this work.
I would like to thank IPMU, Japan for their wonderful hospitality while final part of the work
was done. Most of all I am grateful to the referee for many detailed and helpful comments.
In the review of an early draft, he pointed out that Theorem \ref{t-brd}(b) holds for the new diagram
for $G_{29}$ (suitably modified) and suggested investigating the same question for the other cases
in which algorithm \ref{algo} yields new diagrams. Section \ref{s-newd} is the result of this
investigation.
%
%
%
%
%
\section{The Eisenstein root lattices}
\begin{figure}[ht!]
\begin{center}
\includegraphics{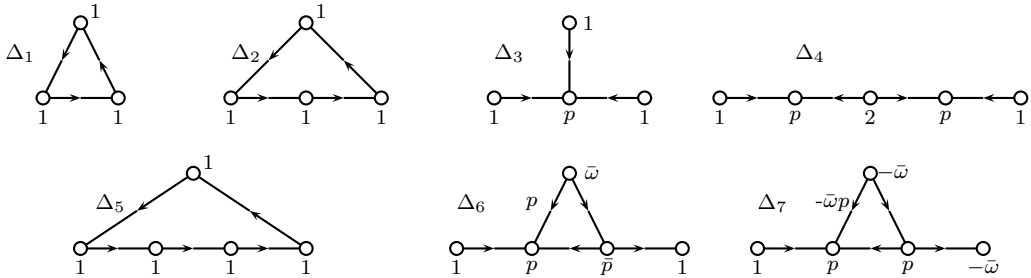}
\caption{A few affine diagrams with balanced numbering.
Choose any vertex $v$ of $\Delta_i$ such that $\Delta_i \setminus \lbrace v \rbrace$
is connected. Then $\Delta_i \setminus \lbrace v \rbrace$ is
a root diagram for $E_4^{\cE}$ if $i = 1$, $E_6^{\cE}$ if $i = 2, 3$ and
$E_8^{\cE}$ if $i = 4, 5, 6, 7$.}
\label{balanced}
\end{center}
\end{figure}

%
%
%
%
\label{seceirootlattice}
Let $\omega = e^{ 2 \pi i/3}$, $\theta = \omega - \bar{\omega}$ and
$p =1  - \bar{\omega}$. Let $\cE = \ZZ[\omega]$.
In this section we shall classify the $\cE$--root lattices, which we define
following Daniel Allcock (see \cite{dja:Y_5552}).
\begin{definition}
An {\it Eisenstein root lattice} or $\cE$--root lattice is a positive definite
$\cE$--lattice $K$, generated by vectors of norm $3$, such that $K \subseteq \theta K'$
(see \cite{dja:Y_5552}).
A root lattice is {\it indecomposable} if it is not a direct sum of two proper
non-zero root lattices.
\par
Let $D$ be a diagram with gram matrix $( \! (m_{i j} ) \!)$.
The following assumptions about $D$ will remain in force for the rest
of this section. We assume that $m_{i j} \in \cE$ for all $i$ and $j$.
We assume that $m_{i i} = 3$ for all $i$. So we omit the labels on the vertices.
If $m_{i j} = -p$, we omit the label on the edge going from $j$ to $i$.
If $m_{i j} = 0$, we omit the edge $\lbrace i,j \rbrace$.
These conventions are adopted when we discuss
connectedness of a diagram.
Each $\cE$--root lattice has at-least one diagram.
Any diagram for an indecomposable $\cE$--root lattice is connected. 
\end{definition}
\begin{remark}
\label{braid-commute}
Let $K$ be a positive definite $\cE$--lattice satisfying $K \subseteq \theta K'$.
The following observations are immediate:
If $r \in K$ has norm $3$, then $r$ is a root of $K$.
The order $3$ reflections in $r$, denoted by
$ \phi_r = \phi_r^{\omega}$ and $\phi_r^{-1} = \phi_r^{\bar{\omega}}$,
belong to the reflection group of $K$.
Let $x$ and $y$ be two linearly independent vectors of $K$ of norm $3$.
The gram matrix of $\lbrace x, y \rbrace$ must be
positive definite, that is, $(\abs{x}^2 \abs{y}^2 - \abs{\ip{x}{y}}^2) > 0$. Since $\ip{x}{y} \in \theta \cE$,
one has $\abs{\ip{x}{y}} \leq \sqrt{3}$.
So either $\ip{x}{y} = 0$, which implies $\phi_x \phi_y = \phi_y \phi_x$, or
$\ip{x}{y} \in  \theta \cE^* $, which implies $\phi_x \phi_y \phi_x = \phi_y \phi_x \phi_y$.
\end{remark}
\begin{definition}
Let $E_{2k}^{\cE}$ be the $\cE$--lattice having a basis $x_1, \dotsb, x_k$ such that
$\abs{x_i}^2 = 3$ for $i = 1, \dotsb, k$,
$\ip{x_i}{x_{i+1}} = -p$ for $i = 1, \dotsb, k-1$ and $\ip{x_i}{x_j} = 0$ if
$j > i+1$. The lattice $E_{2 k}^{ \cE}$ has a diagram
$D_{E_{2k}^{\cE}} = ( x_1 \leftarrow x_2 \leftarrow \dotsb \leftarrow x_k)$.
The $\ZZ$--modules underlying $E_{2}^{\cE}$, $E_{4}^{\cE}$, $E_{6}^{\cE}$ and
$E_{8}^{\cE}$,  with the bilinear form $\tfrac{2}{3}\op{Re}\ip{x}{y}$,
are the lattices $A_2$, $D_4$, $E_6$ and  $E_8$ respectively
(this can be easily checked by computing the discriminant and explicitly exhibiting
the root systems $A_2$, $D_4$, $E_6$ and  $E_8$ inside the real forms of these lattices).
\end{definition}
We use the complex diagrams to give a new proof of the following Theorem (Theorem 2.2 of \cite{dja:Y_5552}):
\begin{theorem}
The only indecomposable Eisenstein root lattices are $E_{2 i}^{\cE}$ with $i = 1, 2, 3, 4$.
\label{thm-classification}
\end{theorem}
The lattices $E_{2}^{\cE}$, $E_{4}^{\cE}$, $E_{6}^{\cE}$ and $E_{8}^{\cE}$
are positive definite. So if some Eisenstein lattice $L$ has the root diagram of
$E_{2j}^{\cE}$, for some $1 \leq j \leq 4$, then $L \simeq E_{2j}^{\cE}$.
Thus it suffices to classify the equivalence classes of root diagrams
of indecomposable $\cE$--root lattices and show that there are only four
classes. The proof of this classification, given below, is like the well-known
classification of $A$-$D$-$E$ root systems.
\begin{definition}
Let $D$ be a connected root diagram for an $\cE$--lattice $L$. 
A {\it balanced numbering} on $D$ is a function from $v(D)$ to
$\cE$, denoted by $x \mapsto n_x$, such that
$n_x = 1$ for some $x \in v(D)$ and
\begin{equation}
\sum_{x \in v(D)} n_x \ip{a}{x} = 0
\label{balancedeq}
\end{equation}
for each $a \in v(D)$.
If $\lbrace n_x \colon x \in v(D) \rbrace$ is a balanced numbering on $D$, then
the vector $y = \sum_{ x \in v(D)} n_x.x \in L$ is orthogonal to each $x \in v(D)$.
So $\abs{y}^2 = 0$. A diagram $D$ is called an {\it affine diagram} if
$D$ admits a balanced numbering but no sub-diagram of $D$ admits one.
\par
Figure \ref{balanced} shows a few affine diagrams, each with a balanced numbering. 
The number shown next to a vertex $x$ is $n_x$. 
Given $a \in v(D)$, suppose $b_1, \dotsb, b_m$ are the vertices connected
to $a$ and $c_j$ is the label on the directed edge going from
$b_j$ to $a$. Then equation \eqref{balanced} becomes
$-3 n_a  = \sum_{j=1}^m  c_j  n_{b_j}$. This is easily verified.
\par
We say that a connected diagram $D$ is {\it indefinite} if $D$
cannot appear as a full sub-graph of a diagram of an Eisenstein
root lattice. Otherwise, we say that $D$ is {\it definite}.
\end{definition}
\begin{lemma}
Let $D$ be a diagram that admits a numbering $\lbrace n_x : x \in v(D) \rbrace$
such that $y = \sum_{ x \in v(D)} n_x .x $ is a norm zero
vector of $L(D)$ and $n_{x_j} = 1$ for some $x_j \in v(D)$. Then $D$ is indefinite.
In particular, if $D$ admits a balanced numbering, then $D$ is indefinite.
\label{indef1}
\end{lemma}
\begin{proof}
Suppose $K$ is an $\cE$--root lattice with a root diagram $D_K$. If
$D \subseteq D_K$, then $y = \sum_{ x \in v(D)} n_x .x  \in K$ has norm $0$.
Since $K$ is positive definite, $y = 0$. If $n_{x_j} = 1$, then
$x_j = - \sum_{ x \neq x_j} n_x .x$, contradicting the fact that $v(D_K)$
is a minimal generating set for $K$. 
\end{proof}
\begin{definition}
Let $D$ be a diagram with vertices $ x_1, \dotsb, x_k $ and
edges $\lbrace x_1 , x_2 \rbrace$, $\lbrace x_2 , x_3 \rbrace, 
\dotsb, \lbrace x_{k-1} , x_k \rbrace$ and $\lbrace x_k , x_1 \rbrace$, that is, a
circuit of length $k$. Changing the vectors $x_1, \dotsb, x_k$ by units if necessary,
we may assume that $\ip{x_i}{ x_{i +1 }} = -p$ for
$1 \leq i \leq k-1$ and $\ip{x_1}{x_k} = -up$ for some $u \in \cE^*$.
We denote this diagram by $\op{Circ}_{k,u}$.
\end{definition}
\begin{lemma}
(a) Let $D = \op{Circ}_{3,-1}$ or $D = \op{Circ}_{3,\bar{\omega}}$. Then $D \simeq
D_{E_6^{\cE}}$.
\par
(b) Let $D = \op{Circ}_{4,\bar{\omega}}$. Then $D \simeq D_{E_8^{\cE}}$.
\par
(c) Suppose $D = \op{Circ}_{k, u}$ with $k = 3, 4$ or $5$ but $D$ is not one of the three
circuits considered in (a) and (b). Then $D$ is indefinite.
\par
(d) Suppose $D$ is one of the diagrams given in figure \ref{balanced}. Then $D$ is indefinite.
\label{indef2}
\end{lemma}
\begin{proof}
(a) Let $v(D)  = \lbrace x_1, x_2, x_3 \rbrace$. Note that $\ip{x_1}{\bar{u} x_3} = \ip{x_1}{x_2} = -p $.
If we take $x_3' = x_2 - \bar{u} x_3$, then $\ip{x_1}{x_3'} = 0$. One checks that
$\abs{x_3'}^2 = 6 + 2 \op{Re}( \bar{u} p ) = 3$ and $ \ip{x_2}{x_3'} = 3 + \bar{u} p \in p \cE^*$,
for $ u = -1$ and $\bar{\omega} $. So $D$ is equivalent to the diagram $D_{E_6^{\cE}}$
formed by the roots $x_1, x_2, x_3'$. 
\par
(b) Let $x_4' = \bar{\omega} x_2 - p x_3 - x_4$. 
One checks that $x_1, x_2, x_3$ and  $x_4'$ form the diagram $D_{E_8^{\cE}}$.
\par
(c)
For $\nu = 1, \dotsb, 6$, let
$(n_{1,\nu}, \dotsb, n_{5,\nu})$ be equal to $(1,-\omega, \bar{\omega}, \bar{\omega},\bar{\omega})$,
$(1, -\omega, \bar{\omega}, -1,-1)$, $(1, -\omega, \bar{\omega}, -1, \omega)$,
$(1, -\omega, \bar{\omega}, -1, \omega)$, $(1,1,1,1,1)$ and $(1, -\omega, -\omega,-\omega,-\omega)$
respectively. One checks that the vector $ y = \sum_{ i = 1}^k n_{i,\nu} x_i$ has
norm zero in $L(\op{Circ}_{k, e^{ 2 \pi i \nu/6}})$, for $k = 3,4, 5$ and $\nu = 1, \dotsb, 6$,
except for the three cases considered in part (a) and (b). 
Part (c) now follows from lemma \ref{indef1}. Part (d) also follows from lemma \ref{indef1}.
\end{proof}
\begin{proof}[Proof of the Theorem \ref{thm-classification}]
Let $D$ be a root diagram for an indecomposable $\cE$--root lattice $K$.
We shall repeatedly use lemma \ref{indef2} in two ways. 
First, it implies that $D$ cannot contain the diagrams
mentioned in part (c) and (d) of the lemma.
Secondly, from the proof of lemma \ref{indef2}, we observe
the following:
\par
{\it If $\op{Circ}_{3, u}$ or $\op{Circ}_{4, u}$ is a sub-graph of $D$, then
we are in one of the cases considered in part (a) or (b) of lemma \ref{indef2} and we can
change one of the vertices to get an equivalent diagram, where one of the
edges has been removed. However this may introduce new edges elsewhere in
the graph.}
\par
Since $K$ is indecomposable, $D$ must be connected.
Let $v(D) = \lbrace x_1 , \dotsb, x_k \rbrace$.
If $k = 1$ (resp. $k=2$), then clearly $K \simeq E_2^{\cE}$ (resp. $K \simeq E_4^{\cE}$).
If $k = 3$, then either $D \simeq D_{E_6^{\cE}}$ or
$D \simeq \op{Circ}_{3,u}$ with $u = -1$ or $\bar{\omega}$, which are
again equivalent to $D_{E_6^{\cE}}$.
\par
Let $k = 4$. We may assume that the diagram $D'$ formed
by $\lbrace x_1, x_2, x_3 \rbrace$ is $D_{E_6^{\cE}}$. If possible,
suppose $D$ is not equivalent to $D_{E_8^{\cE}}$.
Also suppose that $\lbrace x_2, x_4 \rbrace$ is an edge of $D$.
Since $D$ cannot be the affine diagram $\Delta_3$,
either $\lbrace x_1, x_2, x_4 \rbrace$ or $\lbrace x_2, x_3, x_4 \rbrace$
is a circuit in $D$. Without loss, suppose $\lbrace x_1, x_2, x_4 \rbrace$ is
a circuit. Then we can change $x_4$ by adding a multiple of $x_1$ to get
an equivalent diagram where $\lbrace x_2, x_4 \rbrace$ is not an edge.
So $D \simeq E_8^{\cE}$ or $D$ is a circuit. 
In the latter case, $D \simeq \op{Circ}_{4, \bar{\omega}}$, since all other circuits
of length $4$ are indefinite, by lemma \ref{indef2}(c). Lemma \ref{indef2}(b)
implies that $\op{Circ}_{4, \bar{\omega}} \simeq E_8^{\cE}$.
\par
Let $k = 5$. We may assume that the diagram $D'$
formed by $\lbrace x_1, x_2, x_3, x_4 \rbrace$ is $D_{E_8^{\cE}}$.
Since $D$ cannot be the affine diagram $\Delta_4$,
there must be at-least two edges joining $x_5$ with $D'$. Since
the diagrams of the form $\op{Circ}_{5,u}$ are
not definite, $D$ must contain a circuit of length 3 or 4.
If $\lbrace x_1, x_5 \rbrace$ is an edge, it is part of a
circuit of length 3 or 4 and as before, we can remove it by shifting
to an equivalent diagram. The remaining possibilities are
shown in figure \ref{d5cases}. 
\par
\begin{figure}[ht!]
\begin{center}
\includegraphics{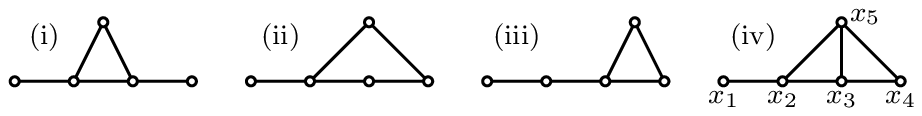}
\caption{An undirected edge between $x_i$ and $x_j$ means $\abs{\ip{x_i}{x_j}} = \abs{p}$.}
\label{d5cases}
\end{center}
\end{figure}
(The arrows on edges are not important here, so they have been omitted).
In cases (iii) and (iv), we may add a multiple of $x_3$ to $x_5$ and disconnect
$x_4$ from $x_5$ (note that this does not introduce an edge between $x_1$ and $x_5$).
So we are reduced to the first two cases. But (i) is affine (either $\Delta_6$ or
$\Delta_7$) and (ii) contains the affine diagram $\Delta_3$ (See figure \ref{balanced}).
\end{proof}
%
%
%
%
%
\section{An attempt to characterize the diagrams for unitary reflection groups}
%
%
%
%
\label{secweylvector}
In this section we want to address a question that was raised in \cite{bmr:cr}:
{\it How to characterize the diagrams for unitary reflection groups?}
We maintain the definitions and notations introduced in \ref{def-reflection},
\ref{def-lattice} and \ref{def-diagram}. 
In \ref{method}-\ref{algo}, we shall describe an algorithm which, given a unitary
reflection group $G$, picks out a set of reflections. As mentioned in the introduction,
these reflections generate $G$ (and
lift to a set of generators for $\op{Braid}(G)$), if the field of definition
of $G$ is $\QQ$ or a quadratic imaginary extension of $\QQ$. There are four fields
to consider, namely, $\QQ(\sqrt{-d})$, where $ d= -1, -2, -3, -7$ and there are $12$ groups
to consider. Most of these can be viewed as sets of automorphisms of certain
complex lattices defined over $\cE$ or $\cG$, namely,
$D_4^{\cE}, E_6^{\cE}, E_8^{\cE}, K_{10}^{\cE}, K_{12}^{\cE}$,
$D_4^{\cG}$ and $E_8^{\cG}$. 
For a description of $E_8^{\cE}$, $E_8^{\cG}$ and $K_{12}^{\cE}$,
see  example 11b, 13b and 10a respectively, in chapter 7, section 8
of \cite{cs:splag}). The lattice $K_{10}^{\cE}$ is the orthogonal
complement of any root in $K_{12}^{\cE}$. 
\par
We take cue from the fact that for a Weyl group, vertices of the Dynkin diagram
correspond to the simple roots, which are the positive roots having
minimal inner product with a Weyl vector.
Our algorithm is based on a generalization of the notion of Weyl vector.
\begin{definition}
\label{def-alpha}
Let $G$ be an irreducible unitary reflection group with a root system $\Phi$ defined over $\cO$.
Let $\Phi_* = \Phi/\cO^*$ be the set of projective roots.
Given a projective root $r$, let $o(r)$ denote the
order of the subgroup of $G$ generated by reflections in $r$.
Let $K$ be the $\cO$--lattice spanned by $\Phi$. 
Let $V$ be the complex vector space underlying $K$ and $\cM = \cup_{ r \in \Phi_*} r^{\bot}$
be the union of the mirrors. Define a function $\alpha: V \setminus \cM \to V$, by
\begin{equation}
\alpha(w) = \sum_{r \in \Phi_*}
o(r)^{-2} \frac{\langle r , w\rangle}{\abs{\langle r , w\rangle} } \frac{r}{\abs{r}}
\label{eq-alpha}
\end{equation}
\end{definition}
\begin{remark}
\label{remark-alpha}
\begin{enumerate}
\item Note that the quantity
$o(r)^{-2} \frac{\langle r , w\rangle}{\abs{\langle r , w\rangle} } \frac{r}{\abs{r}}$
does not change if we change $r$ by a scalar. So the function $\alpha$
is well defined and only depends on the reflection group $G$ and not on the choice
of the roots. The function $\alpha$ descends to a function
$\bar{\alpha} \colon \mathbb{P}(V \setminus \cM) \to \mathbb{P}(V)$.
Also note that $\alpha$ is $G$--equivariant, that is, $\alpha(g w) = g \alpha(w)$
for all $g \in G$. So $\alpha$ induces a function from $\mathbb{P}(V \setminus \cM)/G$
to $\mathbb{P}(V)/G$.
\item  If all the roots have the same norm, then the factor $\abs{r}$
in the denominator is unnecessary and can be omitted from the definition
of $\alpha$.
\item The exponent $-2$ on $o(r)$ was found by experimenting with
the example $G_{26}$, which has reflections of order two and three.
We have included it to indicate one of the ways in which equation
\eqref{eq-alpha} may be modified to possibly include other examples.
If the stabilizer of each mirror has the same order, then
the factor $o(r)^{-2}$ can be omitted from the definition of $\alpha$.
For example, if we consider reflection group of an Euclidean root
lattice (resp. $\cE$--root lattice), then $o(r)$ is always equal
to $2$ (resp. $3$).   
\end{enumerate}
\end{remark}
\begin{topic}{\bf The case of Weyl groups:}
\label{real-case}
Let $W$ be a Weyl group acting on a real vector space $V_{\RR}$ via its standard
representation. Let $\cM_{\RR}$ be the union of the mirrors of $W$.
Then $W$ can be viewed as a unitary reflection group acting on $V = V_{\RR} \otimes \CC$.
\newline
{\it Claim: Let $w, w' \in V_{\RR} \setminus \cM_{\RR}$. 
If $w$ and $w'$ are in the same Weyl chamber, then $\alpha(w) = \alpha(w')$.
Further, $w$ and $\alpha(w)$ belong to the same Weyl chamber.
So $\alpha(w)$ is a fixed point of $\alpha$. }
\begin{proof}
The statement is invariant upon scaling $\alpha$ by a positive factor, so we do
the computation omitting the factor $o(r)^{-2} = 1/4$.
Let $\Phi_+(w)$ be the set of roots having strictly positive inner product
with $w$. Then $\Phi_+(w)$ can be chosen as a set of representatives for
the projective roots, so $\alpha(w) = \sum_{r \in \Phi_+(w)} r/\abs{r}$.
Now, $w$ and $w'$ are in the same Weyl chamber if and only if $\Phi_+(w) = \Phi_+(w')$,
so $\alpha(w) = \alpha(w')$.
\par
To check that $w$ and $\alpha(w)$ belong to the same Weyl chamber, 
first, let $\Phi$ be a root system of type $A_n$, $D_n$, $E_6$, $E_7$ or $E_8$.
Let $w \in V_{\RR} \setminus \cM_{\RR}$.
Then $  \rho =  \tfrac{1}{2} \sum_{ r \in \Phi_+(w)} r$ is a Weyl vector
which belong to the same Weyl chamber as $w$.
Note that $\alpha(w) = \sqrt{2} \rho$ and $\alpha(\sqrt{2}\rho) = \sqrt{2}\rho$.
\par
Now consider the non-simply laced case.
Since $\alpha$ is $G$--equivariant, it is enough to check that $w$ and
$\alpha(w)$ are in the same Weyl chamber, for a single chamber.
We show the calculation for $B_n$. The Weyl group of type $C_n$ is isomorphic to
the group of type $B_n$. Calculations for type $G_2$ and $F_4$ are only little
more complicated and will be omitted.
\par
Let $e_j$ be the $j$--th unit vector in $\ZZ^n$. 
The roots of $B_n$ are $\lbrace \pm e_j , \pm e_i \pm e_j : j < i  \rbrace$.
Choose $w = (w_1, \dotsb, w_n)$ such that $0 < w_1 < \dotsb < w_n$.
Then $ \Phi_+(w) = \lbrace e_j , e_i \pm e_j : j < i  \rbrace$. So
\begin{equation}
\rho^{(n)} = \alpha(w) = 
\sum_{r \colon  \ip{r}{w} > 0 } r/\abs{r} = (1, 1 + \sqrt{2}, 1 + 2 \sqrt{2}, \dotsb, 1 + (n-1) \sqrt{2}).
\label{def-rhok}
\end{equation}
Observe that $\rho^{(n)}$ and $w$ are in the same Weyl chamber. 
So $\alpha(\rho^{(n)}) = \alpha(w) = \rho^{(n)}$.
\end{proof}
\end{topic}
\begin{definition}
Let $\Phi$ be a unitary root system defined over $\cO$. 
In view of \ref{real-case}, a fixed point of $\alpha$ will be called
a {\it Weyl vector}, when the root system is not defined over $\ZZ$.
One may try to find a fixed point of $\alpha$ by iterating the function.
Our method for selecting a set of generating reflections for $G$, is based
on this notion of Weyl vector. Before describing it we show that Weyl
vectors exist.
\end{definition}
%
%
\begin{theorem}
Let $\Phi$ be a root system for a unitary reflection group $G$.
Then the function $\alpha$, defined in \eqref{eq-alpha}, has
a fixed point.
\label{prop-fixpt}
\end{theorem}
For notational simplicity, let $\mu_r = \abs{r}^{-1}. o(r)^{-2}$, so that, 
\begin{equation}
\alpha(w) = \sum_{r \in \Phi_*} \mu_r \abs{\ip{r}{w}}^{-1} \ip{r}{w} r.
\label{alpha-mu}
\end{equation}
The argument given below actually shows, that for any sequence of non-zero positive
numbers $\mu_r$, a function of the form \eqref{alpha-mu} has a fixed point.
We need a lemma, which converts the problem of finding a fixed point of $\alpha$
to a maximization problem.
%
%
\begin{lemma} Let $\Phi$ be a unitary root system and let $\cM$ be the union of mirrors. 
Consider the function $\sym : \PP(V) \to \RR$, defined by
\begin{equation}
\sym (w) = \abs{w}^{-1} \ip{\alpha(w)}{w} = \sum_{r \in \Phi_*} \mu_r \abs{w}^{-1} \abs{\ip{r}{w}}.
\label{sym}
\end{equation}
For $ w \notin \cM$, one has $\partial_w(\sym) = 0$ if and only if
$\alpha(w) = \abs{w}^{-1} \sym(w) w$.
(here $\partial_w$ denotes the holomorphic derivative with respect to $w$).
\label{derivative}
\end{lemma}
\begin{proof} Let $n(w) = \ip{  \alpha(w)}{w}$,
so that $\sym(w) = n(w)/\abs{w} $. We fix a basis for the vector space $V$ and write
$\langle r,w\rangle = r^* M w $, where $r^*$ is the conjugate transpose
of $r$ and $M$ is the matrix of the hermitian form. Differentiating, we get 
$\partial_w( \abs{\langle r,w\rangle}^2) =
\ip{w}{r} r^* M$. So
\begin{equation*}
\partial_w ( \mu_r \abs{\langle r,w\rangle}   )=
 \mu_r ( 2  \abs{\langle r,w\rangle} )^{-1}  \ip{w}{r} r^* M. 
\end{equation*}
Summing over $\Phi_*$, one gets, $\partial_w n(w) = \tfrac{1}{2}\alpha(w)^* M$.
Similarly we get $\partial_w(\abs{w}) =\tfrac{1}{2 \abs{w}} w^* M $.
It follows that 
\begin{equation*}
\partial_w(\sym(w)) = \partial_w\Bigl(  \frac{n(w)}{\abs{w}} \Bigr)  
=\frac{\tfrac{1}{2} \abs{w} \alpha(w)^* M  - \tfrac{1}{2\abs{w} } n(w) w^* M }{\abs{w}^2}
 =\frac{\alpha(w)^* - \frac{\sym(w)}{\abs{w}} w^* }{2\abs{w}} M.
\end{equation*}
Since $M$ is an invertible matrix, the lemma follows.
\end{proof}
%
%
\begin{proof}[Proof of Theorem \ref{prop-fixpt}]
We first prove the following claim:
\newline
{\it Claim: If $w \in V$ lies on a mirror, then the function $\sym$ can not have
a local maximum at $w$.}
\newline
Fix a root $r_0 \in \Phi$ and $w \in r_0^{\bot}$. Assume $\abs{w} = 1$. Take
$w' = w + \epsilon \xi r_0$ , where  $\xi$ is a complex root of unity
and $\epsilon$ is a small positive  real number so that $\epsilon^2$ is negligible. 
We shall show that $\sym(w') > \sym(w)$, for suitable choice of $\xi$.
Ignoring terms of order $\epsilon^2$, we have
$\abs{w'}^2 = \abs{w}^2 + 2\op{Re}(\epsilon \langle w, r_0\rangle) = \abs{w}^2 = 1$.
Let $\Psi_0 = \Phi_* \cap w^{\bot}$ and $\Psi_1 = \Phi_* \setminus \Psi_0$.
Then
\begin{equation*}
\sym(w')
= \sum_{r \in \Psi_0 } \mu_r \abs{\langle r,w'\rangle}
+ \sum_{r \in \Psi_1 } \mu_r \abs{\langle r,w'\rangle} 
= \epsilon \sum_{\Psi_0}  \mu_r \abs{\langle r,r_0\rangle}
+  \sum_{ \Psi_1 } \mu_r \bigl| \langle r,w\rangle +\epsilon \xi \ip{r}{r_0} \bigr|.
\end{equation*}
Let $a = \ip{r}{w} \neq 0$ and $b = \xi\ip{r}{r_0} $.
Using the first order expansion of $(1 + x)^{1/2}$, we have,
\begin{equation*}
\abs{ a + \epsilon b} = \sqrt{\abs{a}^2 + 2 \op{Re}( \epsilon b \bar{a} )}
= \abs{a} +  \op{Re}( \epsilon b \bar{a}\abs{a}^{-1} ).
\end{equation*}
It follows that
\begin{equation*}
\sym(w') = \sym(w) + \epsilon \sum_{\Psi_0} \mu_r \abs{\langle r,r_0\rangle}
+\epsilon \op{Re} \Bigl( \xi  \sum_{\Psi_1} 
 \mu_r \ip{r}{r_0} \overline{\langle r, w\rangle} \abs{\langle r,w\rangle}^{-1} \Bigr).
\end{equation*}
Note that the third term can be made non-negative by choosing $\xi$ suitably, and
the second term is positive, since $r_0 \in \Psi_0$.
This proves the claim.
\par
The function $\sym$ is continuous on $\PP(V)$, so it attains its global maximum,
say at $w_0$. The claim we just proved implies that
$w_0 \notin \cM $, so $\partial_{w_0}(\sym) = 0$. Let $w_1 =
 \abs{w_0}^{-1} \sym(w_0) w_0$. Lemma \ref{derivative} implies that
$\alpha(w_1) = \alpha(w_0) = w_1 $. 
\end{proof}
%
%
\begin{definition} 
Let $(\cO, G, \Phi, K, V)$ be as in \ref{def-alpha}. Given $w \in V$, let
$(r_1, \dotsb, r_N)$ be the projective roots of $G$ arranged so that
$d(r_1^{\bot}, w) \leq \dotsb \leq d(r_N^{\bot}, w)$, where $d$ is the
Fubini-Study metric on $\PP(V)$. So
\begin{equation*}
d(r^{\bot}, w) = \sin^{-1} ( \abs{\ip{r}{w}}/\abs{r} \abs{w}).
\end{equation*}
Let $k$ be the minimum number of reflections needed to generate $G$. 
Define $\Delta(w) = \lbrace r_1, \dotsb, r_k \rbrace$. In other words,
$\Delta(w)$ consists of $k$ projective roots, whose mirrors are closest to $w$.
\end{definition}
\begin{topic}{\bf Method to obtain a set of ``simple reflections'': }
\label{method}
We now describe a computational procedure in which, the input is a
unitary reflection group $G$, (or equivalently, a projective root system $\Phi_*$ for $G$
and the numbers $\lbrace o(r) \colon r \in \Phi_* \rbrace$)
and the output is either the empty set or a non-empty set of projective roots,
to be called the {\it simple roots}. 
The reflections in the simple roots are called {\it simple reflections}.
A set of simple reflections form a {\it simple system}.
\par
{\it 
\begin{enumerate}
\item From the set $\lbrace w \colon \alpha(w) = w \rbrace$, choose $w$ such that $\sym(w)$,
defined in \eqref{sym}, is maximum. 
\item Let $\Delta(w) = \lbrace r_1, r_2, \dotsb, r_k \rbrace$.
\item If $r_1, \dotsb, r_k$ are linearly dependent,
then return the empty set.\footnote{If the group $G$ is not well generated,
then one should use obvious modifications; e.g, step (3) should be rephrased as follows:
If $\op{dim}(\op{span}\lbrace r_1, \dotsb, r_k \rbrace)$ is less than the rank of 
$G$, then return the empty set.}
\item If $r_1, \dotsb, r_k$ are linearly independent, then return 
$\Delta(w)$ (the simple roots).
\end{enumerate}
} 
In practice, for each group $G$ to be studied, we execute the following algorithm many times.
\end{topic}
\begin{topic}{\bf Algorithm:}
\label{algo}
{\it Start with a random vector $w_0 \in V$. Generate a sequence $w_n$ by $w_n = \alpha(w_{n-1})$.
If the sequence stabilizes, then say that the algorithm converges and let $w  = \lim w_n$.
Note down $\Delta(w)$ and $\sym(w)$.}
\newline
\newline
For computer calculation, we assume that $w_n$ stabilizes if $\abs{w_{n+1} - w_n}^2/\abs{w_n}^2$ 
becomes small, say less than $10^{-8}$, and remains small and decreasing for many successive
values of $n$. Note that if $\alpha(w) = w$, then $\sym(w) = \abs{w}$. From the values
of $\sym(w)$, the maximum, denoted by $\sym_{\mathrm{max}}^G$, can be found. 
(In all the examples that we have studied, $\sym(w)$ takes at-most two values on the
set of fixed points of $\alpha$ found experimentally, so finding
the maximum is not difficult).
Each instance of the algorithm, that produced a vector $w$ with $\sym(w) = \sym_{\mathrm{max}}^G$,
now yields a simple system $\Delta(w)$ provided that $\Delta(w)$ is a linearly
independent set.
\end{topic}
\begin{topic}{\bf Observations for Weyl groups: }
First, suppose that $\Phi$ is a root system for a Weyl group $W$. We maintain the notations
of \ref{real-case} and ignore the factor $o(r)^{-2} = \tfrac{1}{4}$
in the definition of $\alpha$.
Given $w \in V_{\RR}$, let $\delta_w(r) = \abs{\ip{r}{w}}/\abs{r}$.
Arranging the roots according to increasing distance from $w$ 
(according to the spherical metric or the Fubini-Study metric) is equivalent
to arranging them according to increasing order of $\delta_w(r)$.
\newline
{\it Claim: Let $W$ be a Weyl group, $w_0 \in V_{\RR} \setminus \cM_{\RR}$
and $w = \alpha(w_0)$. 
Then the algorithm \ref{algo} converges in one iteration and yields a simple
system $\Delta(w)$. For Weyl groups, the definition of a set of simple
roots given in \ref{method} agrees with the classical notion.
Further, for each simple root $r \in \Delta(w)$, we have $\delta_w(r) = 1$.
So the simple mirrors are equidistant from $w$.}
\begin{proof}
We saw in \ref{real-case} that $w =\alpha(w_0)$ is a fixed point of $\alpha$,
so $w = \alpha(w_0) = \alpha^2(w_0) = \dotsb$, that is, the algorithm converges in
one iteration. 
\par
Suppose $\Phi$ is a root system of type $A$, $D$ or $E$.
Then $w = \alpha(w_0) = \sqrt{2} \rho$, where $\rho$ is a Weyl vector.
So $\Delta(w) = \Delta(\rho)$ consists of a set of simple roots and 
for each $r \in \Delta(w)$, one has $\ip{r}{\rho} = 1$, so 
$\delta_w(r) = 1$.
\par
In type $B_n$, if $w_0$ is in the Weyl chamber containing the vector $\rho^{(n)}$ given
in equation \eqref{def-rhok}, then one has $ w = \rho^{(n)} = \alpha(\rho^{(n)}) = \dotsb $. 
Observe that the function $\delta_w(r)$ attains its minimum for the simple roots
$\Delta(w) = \lbrace e_1, e_2 - e_1, \dotsb, e_n - e_{n-1} \rbrace$
and $\delta(r) = 1$ for each $r \in \Delta(w)$. 
The claim was verified for $G_2$ and $F_4$ without difficulty.
\end{proof}
\end{topic}
\begin{topic}{\bf Observations for complex root systems:}
In the following discussion, let $G$ be one of the unitary reflection groups
from the set $\lbrace G_{4}, G_{5}, G_{25}, G_{26}, G_{32}, G_{33}, G_{34} \rbrace
\cup \lbrace G_8, G_{29}, G_{31}\rbrace \cup \lbrace G_{12} \rbrace \cup \lbrace G_{24} \rbrace$.
The groups given in four subsets are defined over
$\QQ(\sqrt{-3})$, $\QQ(\sqrt{-1})$, $\QQ( \sqrt{-2})$ and $\QQ(\sqrt{-7})$
respectively. In each case, a projective root system $\Phi_*(G)$ for $G$ is chosen.
These are described in the appendix \ref{ap-1}.
For each of these groups $G$, we have run algorithm \ref{algo} at-least one thousand
times and obtained many simple systems by method \ref{method}.
The calculations were performed using the GP/PARI calculator. 
The main observations made from the computer experiments are the following:
\newline
\newline
{\it  The sequence $\lbrace w_n \rbrace $ stabilizes in each trial for each
$G$ mentioned above. 
The simple reflections generate $G$ and satisfy the same set of relations $R_G$, every time.
In other words, the simple roots form the same ``diagram" every time.
For $G_4$, $G_5$, $G_8$, $G_{25}$, $G_{26}$, $G_{32}$, these are Coxeter's diagrams, (see \cite{bmr:cr}).
For $G_{24}$, the relations $R_G$ are given in \ref{a-5} (this is the third presentation given in \cite{bm:ep}).
For $G_{29}$, $G_{31}$, $G_{33}$ and $G_{34}$ new
diagrams are obtained.  These diagrams are given in figure \ref{fig-3} and corresponding
presentations for the braid groups are given in \ref{t-rel}, \ref{p-mutate}. Finally, for
$G_{12}$\footnote{In the notation of section \ref{s-newd}
 the Coxeter diagram for $G_{12}$ is a triangle with each edge marked with $\infty$ and for $G_{24}$
it is a triangle with two double edges and one single edge. We have not drawn these.
We should remark that in these two examples, we could not find presentations, on our generators, 
consisting of only cyclic homogeneous relations, though such presentations exist (see \cite{bmr:cr}, \cite{bm:ep}). }
a presentation is given in \ref{a-4}.  
}
\newline
\newline
Further observations from the computer experiments are summarized below.
\begin{enumerate}
\item 
If $G$ is not $G_{29}$, $G_{31}$ or $G_{32}$, then for all $w$ such that $\alpha(w) = w$,
the function $S(w)$ attains the same value. 
So each  trial of the algorithm yields a maxima $w$ for $\sym$.
For $G_{29}$, $G_{31}$ and $G_{32}$,
the function $\sym(w)$ attains two values on the fixed point set of $\alpha$.
In these three cases, the $\ZZ$--span of the roots form the $E_8$ lattice.
%
\item If $G \neq G_{33}$ and $G$ is well generated, then in each trial of the algorithm,
we find that the vectors $\lbrace r_1, \dotsb, r_k \rbrace$ are linearly
independent. So each choice of a maxima $w$ for the function $\sym$ yields
a set of simple roots $\Delta(w)$.
For $G_{33} = R(K_{10}^{\cE})$, in most of the trials,
we find that  $\lbrace r_1, \dotsb, r_5 \rbrace$ are linearly dependent.
So most trials do not yield a set of simple roots $\Delta(w)$.
In an experiment with $5000$ trials, only $297$ yielded simple systems.
\item 
Let $G\in \lbrace G_{33}, G_{34}, G_{29} \rbrace$.
These are the three well generated groups for which the diagrams obtained by method \ref{method}
are different from the ones in the literature. 
Let $g_1, \dotsb, g_k$ be a set of simple reflections of $G$ obtained by method
\ref{method}. For these three groups, we have verified the following result.
(almost a re-statement of \ref{t-brd}): 
\newline
\newline
{\it Fix a permutation $\pi$ so that the order of the product
$p = g_{\pi_1} \dotsb g_{\pi_k}$ is maximum (over all permutations).
Then, the order of $p$ is equal to the Coxeter number of $G$ (denoted by $h$) and
the eigenvalues of either $p$ or $p^{-1}$ are
$e^{ 2 \pi i (d_1 - 1)/h }, \dotsb, e^{ 2 \pi i (d_k - 1)/h }$,
where $d_1, \dotsb, d_k$ are the invariant degrees of $G$.
}
\newline
\newline
The invariant degrees of $G_{29}$, $G_{33}$ and $G_{34}$ are $(4,8,12,20)$,
$(4,6,10,12,18)$ and $(6,12,18,24,30,42)$ respectively. In all three cases,
$h$ is equal to the maximum degree. 
\item Start with $w_0 \in V$ and consider the sequence defined
by $w_n = \alpha(w_{n-1})$. Roughly speaking, each iteration of the function $\alpha $
makes the vector $w_n$ more symmetric with respect to the set of projective roots $\Phi_*$.
The function $\sym$ measures this symmetry.
So the fixed points $w$ of $\alpha$ are often the vectors that are
most symmetrically located with respect to $\Phi_*$. 
\par
Based on the above discussion, an alternative definition of a Weyl
vector may be suggested, namely, a vector $w \in V$, such that $d(w, \cM)$ is maximum.
(This was suggested to me by Daniel Allcock).
It seems harder to compute these vectors, so we have not
experimented much with this alternative definition.
However, we would like to remark that for some complex and quaternionic
Lorentzian lattices, similar analogs of Weyl vectors and simple roots, are useful.
(One such example is studied in \cite{tb:el}; other examples are studied in
\cite{tb:thesis}. In these examples of complex and quaternionic Lorentzian lattices,
the simple roots are again defined as those whose mirrors are closest to the ``Weyl vector".) 
\item 
The method \ref{method} fails for the primitive unitary reflection groups
that are not defined over $\QQ$ or a imaginary quadratic extension of $\QQ$
and for the imprimitive groups $G(de, e, n)$, except
when they are defined over $\QQ$, that is, 
for the cases $G(1,1,n+1) \simeq A_n$, $G(2,1,n) \simeq BC_n$ and
$G(2,2,n) \simeq D_n$.
It fails in the sense that the set $\Delta(w)$ does not in general
form a minimal set of generators for the group. This was found by
experimenting with $G(de,e,n)$ for small values of $(d,e,n)$ and
also with $G_6$ and $G_9$ (see appendix \ref{a-4}).
Although method \ref{method} fails, the following observation holds
for $G(de,e,n)$:
\par
{\it
Let $k$ be the minimum number of reflections needed to generate $G(de,e, n)$
($k = n$ or $n+1$).
There exists a vector $w \in V$ such that, if $r_1^{\bot}, \dotsb, r_k^{\bot}$
are the mirrors closest to  $w$, then reflections in
$\lbrace r_1, \dotsb, r_k \rbrace$ generate $G(de, e, n)$.
}
\par
For $G(de,e,n)$, one can take $w = \rho^{(n)}$ (the vector we obtained for the Weyl group $B_n$;
see equation \eqref{def-rhok}).
It is easy to check that $\Delta(\rho^{(n)})$ forms the known diagrams
for $G(d e, e, n)$.
We found $\rho^{(n)}$ by using an algorithm that tries to find a point in $V$
whose distance from $\cM$ is at a local maximum.
For small values of $d$, $e$ and $n$, we find that this algorithm always
converge to  $\rho^{(n)}$. 
%
\end{enumerate}
\end{topic}
%
%
%
%
%
\section{The braid groups for $G_{29}$, $G_{31}$, $G_{33}$, $G_{34}$.}
\label{s-newd}
\begin{figure}
\begin{center}
\includegraphics{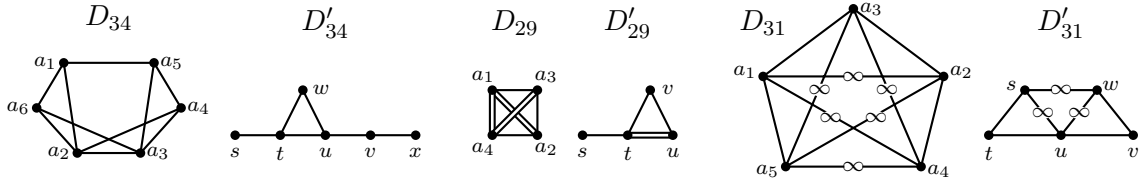}
 \caption{ Let $D_{33} = D_{34} \setminus \lbrace a_6 \rbrace$ and $D'_{33} = D'_{34} \setminus \lbrace x \rbrace$.
 The diagrams $D_{N}$ and $D'_{N}$  encode presentations for $\op{Braid}(G_N)$.}
\label{fig-3}
\end{center}
\end{figure}

%
%
%
%
\begin{topic}
For this section, let $N \in \lbrace 29, 31, 33, 34 \rbrace$.
Algorithm \ref{algo}, applied to $G_N$, yields new diagrams, denoted by $D_N$ (see figure \ref{fig-3}).
Before stating the main result, proposition \ref{p-mutate}, we need some notations.
\end{topic}
\begin{topic}{\bf Notations:}
Let  $\lbrace x_i \colon  i \in \ZZ/k \ZZ \rbrace$ be elements of a monoid $M$.
Let $\rel{m}{x_0, \dotsb, x_{k-1} }$ denote the positive and homogeneous relation
\begin{equation} 
x_0 x_1 \dotsb x_{m-1} = x_1 x_2 \dotsb x_{m}.
\label{e-def}
 \end{equation}
 For example $\rel{2}{a,b}$ (resp. $\rel{3}{a,b}$) says that $a$ and $b$ commutes (resp. braids);
while $\rel{5}{a, b}$ (resp. $\rel{4}{a,b,c}$) stands for the relation $a b a b a = b a b a b$ and
(resp. $a b c a  = b c a b $). 
If $x, y$ are elements in a group, let $c_x(y) : = x y x^{-1}$.
\par
Given a diagram $D$ (such as in figure \ref{fig-3}) let $\op{Cox}(D, \infty)$ be the group defined by 
generators and relations as follows: The generators of $\op{Cox}(D, \infty)$ correspond to the vertices of $D$.
The relations are encoded by the edges of $D$:  $k$ edges between vertices $p$ and $q$ encodes
 the relation $\rel{k+2}{p,q}$. In particular, no edge between $p$ and $q$ indicates that $p$ and $q$
 commute, while an edge marked with $\infty$ indicates that $\op{Cox}(D, \infty)$ has
 no defining relation involving only the generators $p$ and $q$.
Let $\op{Cox}(D, n)$ be the quotient of $\op{Cox}(D, \infty)$ obtained by imposing the relation
$p^n = 1$ for each vertex $p$ of $D$. Let $\op{Cox.Rel}(D)$ be the defining relations of $\op{Cox}(D, \infty)$.
\par
Let $\tilde{A}_n$ denote the affine Dynkin diagram of type $A_n$. (Picture it as a regular polygon with
$(n+1)$ vertices).  
Fix an automorphism $\rho$ of $\tilde{A}_n$ (hence of $\op{Cox}(\tilde{A}_n, \infty)$),
that rotates the Dynkin diagram $\tilde{A}_n$ by an angle $2 \pi/n$. Let $x$ be a vertex of $\tilde{A}_n$.
We say that the
relation $\rel{m}{x, \rho(x), \dotsb, \rho^{n}(x)}$ is {\it cyclic} in $\op{Cox}(\tilde{A}_n, \infty)$ if
the relation $\rel{m}{\rho(x), \dotsb, \rho^{n}(x), x}$ holds in the quotient
$\op{Cox}(\tilde{A}_n, \infty)/ \langle \rel{m}{x, \rho(x), \dotsb, \rho^{n}(x)} \rangle $
(so $\rel{m}{y, \rho(y), \dotsb, \rho^{n}(y)}$ holds
for each vertex $y \in \tilde{A}_n$).
\par
The two lemmas stated below help us verify proposition \ref{p-mutate}. But these
might be of independent interest for studying groups satisfying relations of the form \eqref{e-def}.
The proofs are straight-forward and given in \ref{pf-l-cyclic} and \ref{pf-l-rel}. It is easy to write
down more general statements than those stated below and give an uniform proof, at-least for Lemma \ref{l-rel}.
To keep things simple, we have resisted this impulse to generalize.
\end{topic}
\begin{lemma} Let $(x_0, \dotsb, x_n) = (x, \rho(x), \dotsb, \rho^{n}(x))$ be the vertices of $\tilde{A}_n$.
\par
 (a) $\rel{n}{x_0, \dotsb, x_n }$ is cyclic in $\op{Cox}(\tilde{A}_n, \infty)$
 for all $n$.
 \par
  (b) $\rel{n+2}{ x_0, \dotsb, x_n}$ (resp. $\rel{2 n+2}{ x_0, \dotsb, x_n}$) 
  is cyclic in $\op{Cox}(\tilde{A}_n, \infty)$ if and only if $x_n = x_2$ (resp. $x_{n-1} = x_1$). In particular
 $\rel{4}{x_0, x_1, x_2}$ and $\rel{6}{x_0, x_1, x_2}$ are cyclic in $\op{Cox}(\tilde{A}_2, \infty)$. 
 \par
 (c) Let $n > 2$. Then $\rel{2 n +3}{ x_0, \dotsb, x_n}$ is cyclic in $\op{Cox}(\tilde{A}_n, \infty)$ if and only if $x_{n-1} = x_2$. In particular
 $\rel{9}{x_0, x_1, x_2, x_3}$ is cyclic in $\op{Cox}(\tilde{A}_3, \infty)$. 
 \label{l-cyclic}
\end{lemma}
\begin{lemma} (a) Assume $x$ braids with $y$ in a group $G$. Then the following are equivalent:
\par
(i) $\rel{4}{x, y, z}$ (ii) $c_{x}(y)$ commutes with $z$. (iii) $c_y(z)$ commutes with $x$.
\par
If $x$ also braids with $z$ then (i),(ii), (iii) are equivalent to (iv): $c_z(x)$ commutes with $y$.
\par
(b) Assume $x$ braids with $y$ and $z$ in a group $G$. Then the following are equivalent:
\par
(i) $\rel{6}{x,y,z}$ (ii) $c_x(y)$ braids with $z$. (iii) $c_y(z)$ braids with $x$. (iv) $c_z(x)$ braids with $y$.
\par
(c)  Suppose $x, y, z, w$ are the Coxeter generators of $\op{Cox}(\tilde{A}_3, \infty)$. Then $\rel{9}{ x, y, z ,w }$
holds if and only if $c_{x y }(z)$ braids with $w$.
\label{l-rel}
\end{lemma}
\begin{topic}{\bf The relations: } 
Let $N \in \lbrace 29, 31, 33, 34 \rbrace$.
The vertices of the diagram $D_N$ and $D'_N$ given in figure \ref{fig-3} correspond to generators of $\op{Braid}(G_N)$.
The edges indicate the Coxeter relations. However some more relations are needed to obtain a presentation of $\op{Braid}(G_N)$.
These relations are given below. All of them are of the form \eqref{e-def}.
\begin{align*} 
 E_{34}& 
 = \lbrace \rel{4}{a_2, a_3, a_4}, \rel{4}{a_3,a_4,a_5}, \rel{9}{a_3, a_2,a_1,a_5}, \\
 & \;\;\;\; \;\;\;\; \rel{4}{a_2, a_1, a_6}, \rel{4}{a_2, a_3, a_6}, \rel{9}{a_1, a_5, a_3, a_6} \rbrace, \;\;\; \;\;\;\;\;\;
 E'_{34} = \lbrace \rel{6}{t, u, w}\rbrace,\\
 E_{29} &= \lbrace \rel{4}{ a_4, a_2, a_1},  \rel{4}{ a_3, a_1, a_4},  \rel{6}{ a_3, a_1, a_2} \rbrace, \;\;\;\;\;\;\;\; \;\;\; \;\;\;\;\;
 E'_{29} =\lbrace \rel{6}{u,t,v} \rbrace, \\
 E_{31} &=  \lbrace \rel{4}{a_4, a_2, a_1}, \;\;  \rel{6}{ a_2, a_3, a_4},   \rel{3}{a_3, a_5,a_4},   \\
  & \rel{3}{a_5, a_4,a_3}, \rel{4}{a_1, a_5, a_2},   \rel{4}{a_3, a_1, a_4}, \rel{4}{a_4, a_2, a_5}\rbrace,\;\;
 E'_{31} = \lbrace \rel{3}{s, u, w}, \rel{3}{u,w,s} \rbrace.
  \end{align*}
 Let $E_{33}$ (resp. $E'_{33})$ be the relations in $E_{34}$ (resp. $E'_{34}$) that does not involve $a_6$ (resp. $x$).
 Let $R_N  = \op{Cox.Rel}(D_N) \cup E_N$ and  $R'_N = \op{Cox.Rel}(D'_N) \cup E'_N$.
Let $B_N$ (resp. $B'_N$) be the group generated by the vertices of the diagram $D_N$ (resp. $D'_N$) satisfying 
the relations $R_N$ (resp. $R'_N$). It was conjectured in \cite{bm:ep} and proved in \cite{db:fc} that
$\op{Braid}(G_N)  \simeq B'_N$.
\label{t-rel}
\end{topic}
 \begin{proposition}Assume the setup given in \ref{t-rel}. 
The presentations $B_N$ and $B'_N$ are equivalent. So $B_N$ gives a presentation of $\op{Braid}(G_N)$. 
The quotient of $B_N$ obtained by imposing the relations $a_i^2 = 1$, for all $i$, is isomorphic to $G_N$.
 \label{p-mutate}
\end{proposition}
\begin{proof}[sketch of proof]
We define the maps  $\varphi_N: B'_N \to B_N$ and $\psi_N: B_N \to B'_N$ on the generators. Let
\begin{equation*} 
\varphi_{29}: (s, t, v, u) \mapsto (c_{(a_2 a_3 a_1)^{-1}}(a_4),  a_1,  a_2,  a_3), \;\;\;\; \psi_{29} : (a_1, a_2, a_3, a_4) \mapsto ( t, v, u, c_{ v u t} (s) ).
 \end{equation*}
\begin{align*} 
& \varphi_{31}: (s, t, u, w, v) \mapsto (a_5, a_1, a_3, c_{a_3^{-1}}(a_4),  c_{a_4}(a_2) ),\\
& \psi_{31} : ( a_5, a_1, a_3, a_4, a_2) \mapsto ( s,  t, u, c_u(w) ,  c_{ u w^{-1} u^{-1} }(v)).
 \end{align*}
 \begin{align*} 
& \varphi_{34}: (s, t, u, v, w) \mapsto (a_1, a_2, a_3, c_{a_3}(a_4), c_{a_2 a_1} (a_5), c_{a_4 a_2 a_1 a_5 a_3}(a_6)), \\
& \psi_{34}: (a_1, a_2, a_3, a_4, a_5) \mapsto (s, t, u, c_{u^{-1}}(v), c_{(ts)^{-1} }(w), c_{ (v u w t s u )^{-1}}(x)).
\end{align*}
Let $\varphi_{33}$ and $\psi_{33}$ be the restrictions of $\varphi_{34}$ and $\psi_{34}$ respectively.
 \par
To check that $\varphi_N$ (resp. $\psi_N$) is a well defined group homomorphism, we have to verify that
$(\varphi_N(s), \varphi_N(t), \dotsb, )$ (resp. $(\psi_N(a_1), \psi_N(a_2), \dotsb)$) satisfy the relations $R'_N$ (resp. $R_N$).
This verification was done by hand using  lemmas \ref{l-cyclic} and \ref{l-rel}.
The details, given in \ref{pf-g33}, \ref{pf-g29} and \ref{pf-g31}, are rather tedious.
It is easy to see that $\varphi_N$ and $\psi_N$ are mutual inverses, so $B_N \simeq B'_N$. 
\par
The generators $a_i$ for $B_N$ were found as follows.
We first found  generators $\bar{a}_i$ of order $2$ in $G_N$ using
algorithm  \ref{algo} and then let $a_i$ be a lift of $\bar{a}_i$, that is, we chose 
$R_N$ to be an appropriate subset of the relations satisfied by $\lbrace \bar{a}_1, \bar{a}_2, \dotsb \rbrace$ in $G_N$.
So, from our construction, we know that there are reflections of order $2$ in $G_N$ that satisfy the relations $R_N$, that
is, $G_N$ is a quotient of $\bar{B}_N := B_N/\langle a_i^2 = 1 \text{\; for all \;} i \rangle$. Using coset enumeration on
the computer algebra system MAGMA, we verified that $\bar{B}_N$ has the same order as $G_N$. 
\end{proof}
\begin{remark} Here are some concluding remarks for this section.
\begin{enumerate}
\item The transformations, $\varphi_N$ and $\psi_N$, given in the proof of \ref{p-mutate} were found first by computing inside
the reflection group $G_N$ and then choosing appropriate lifts to $\op{Braid}(G)$.  
\item It is interesting to note that all the relations encountered in the presentations of the braid groups of type $G_{33}$
and $G_{34}$ are cyclic relations (see \ref{l-cyclic}) of the form $\rel{n}{x_1, \dotsb, x_m}$
where $x_1, \dotsb, x_m$ is a minimal cycle (actually a triangle or a square) inside the diagrams. 
\item  The relation $\rel{9}{a_4, a_2, a_1, a_5}$ (resp. $\rel{6}{ a_2, a_3, a_4}$) hold as a consequence of the relations $R_{33}$ (resp. $R_{29}$)
(see \ref{r-remainingrel}).
So the content of Prop. \ref{p-mutate} for $N = 29, 33$ and $34$ can be succinctly stated as follows:
The relations needed to present $\op{Braid}(G_N)$ are $\op{Cox.Rel}(D_N)$ and one relation of the
form $\rel{n}{x_1, \dotsb, x_m}$ for each minimal cycle $x_1, \dotsb, x_m$ in $D_N$.  Further, this $n$ is the smallest
integer for which the relation $\rel{n}{x_1, \dotsb, x_m}$ holds in $G_N$.
It might be interesting to find out all the unitary reflection groups for which a similar statement is true.
\item The situation with the non-well generated group $G_{31}$ is not as nice. First of all, we were not able to verify
the previous remark for $\op{Braid}(G_{31})$. Further, deleting the relations involving $a_5$ from $R_{31}$
only yields a proper subset of $R_{29}$.
The problem is that the relations $\rel{4}{a_1, a_2}, \rel{4}{a_1,a_4},
\rel{4}{a_2,a_5}, \rel{4}{a_3, a_4}, \rel{4}{a_3, a_5}, \rel{4}{a_4, a_5}$ hold in $G_{31}$, but
we were not able to check whether these relations are implied by $R_{31}$. This would be equivalent to checking
whether the relations $\rel{4}{s,u}, \rel{4}{u,w},\rel{4}{s,w}$ are implied by $R'_{31}$.
\item Some relations of the form  $\rel{k}{x_0, \dotsb, x_n}$ were studied in \cite{ccs:26}.
Conway called them deflation relations, because they often ``deflate" infinite Coxeter groups to finite groups. 
The groups studied in this section provide some examples of this. We mention two other examples: 
(i) the quotient of the affine Weyl group $\op{Cox}(\tilde{A}_n, 2)$ obtained by the 
adding the cyclic (see \ref{l-cyclic}(a)) relation $\rel{n}{x_0, \dotsb, x_n}$ is
the finite Weyl group $\op{Cox}(A_n, 2)$. 
(ii) There is a graph $D$ with $26$ vertices such that the quotient of
the infinite Coxeter group $\op{Cox}(D,2)$ obtained by adding the cyclic deflation relations
$\rel{11}{x_0, \dotsb, x_{10}}$, for each minimal cycle with $12$ vertices in $D$, is the wreath product
of the monster simple group with $\ZZ/2 \ZZ$ (see \cite{ccs:26}).
\end{enumerate}
\end{remark}
%
%
%
%
%
%
%
\section{the affine reflection groups}
\begin{figure}
\begin{center}
\includegraphics{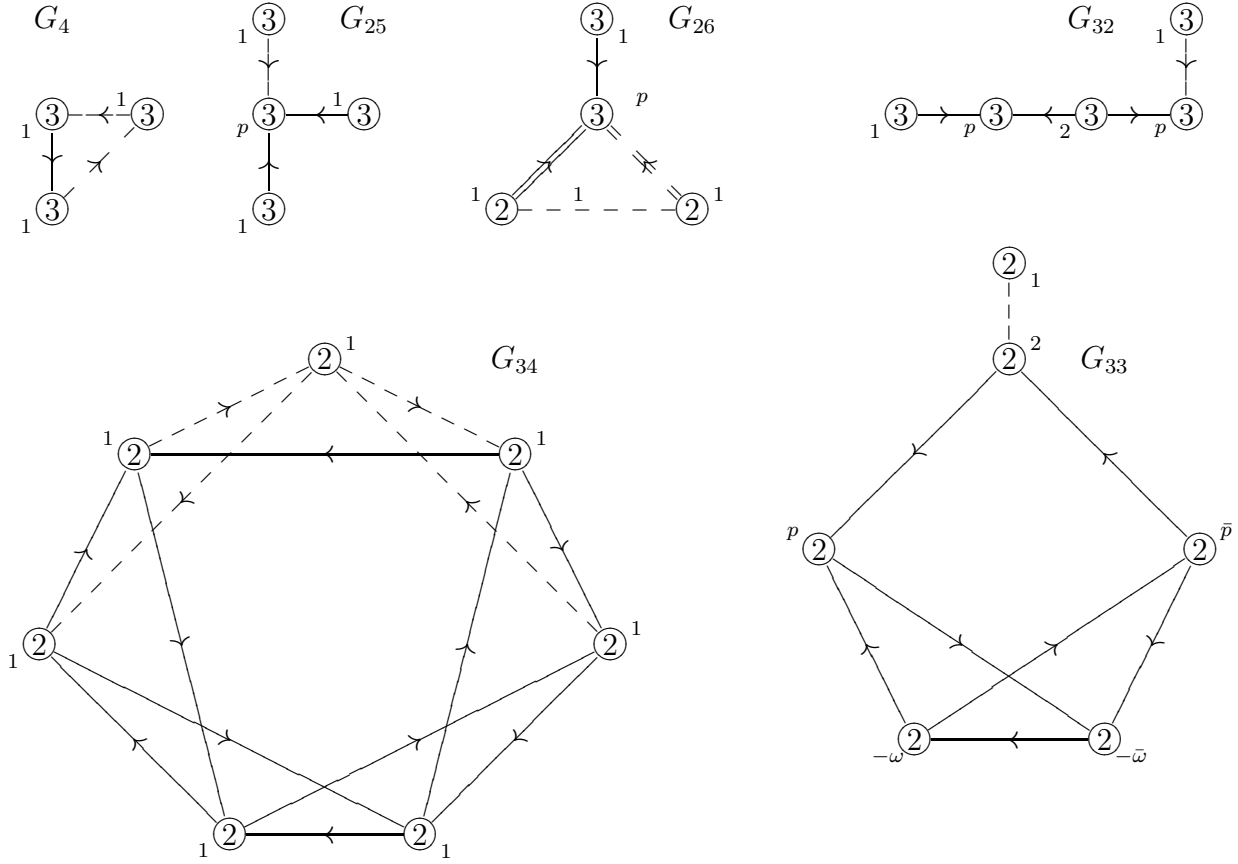}
\caption{Affine diagrams with balanced numbering (shown next to the vertices).
The number shown inside the vertex is the norm of the root as well as
the order of a reflection in that root. 
An edge (resp. a double edge) between $x$ and $y$ implies the
Coxeter relation $\phi_x \phi_y \phi_x = \phi_y \phi_x \phi_y$ (resp.
$\phi_x \phi_y \phi_x \phi_y = \phi_y \phi_x \phi_y \phi_x$).
An unmarked directed edge or double edge from $y$ to $x$ means that
$\ip{x}{y} = -p$ for $G_4$, $G_{25}$, $G_{26}$, $G_{32}$ and 
it means $\ip{x}{y} = \omega $ for $G_{33}$, $G_{34}$.}
\label{fig-affine}
\end{center}
\end{figure}

%
%
%
%
\label{secaffine}
\begin{topic}
In this section we shall describe affine diagrams for the 
primitive unitary reflection groups defined over $\cE$, except for $G_5$.
(Including $G_5$ would further complicate notations).
An affine diagram is obtained by adding an extra node to
the corresponding ``unitary diagram". Each affine diagram admits a
balanced numbering. 
A unitary diagram can be extended to an affine diagram in
many ways. We have chosen one that makes the diagram more
symmetric.  The Weyl vector, that yielded the unitary diagram,
is often fixed by the affine diagram automorphisms.
\par
The discussion below and the lemma following it are direct
analogs of the corresponding results for Euclidean root lattices. We have
included a proof since we could not find a convenient reference. 
\par
Let $\Phi$ be a unitary root system defined over $\cE$, for an unitary reflection group $G$.
Let $K$ be the $\cE$--lattice spanned by $\Phi$. Assume that the subgroup of $\Aut(K)$ generated by
reflections in $\Phi$ is equal to $G$.
Let $\cE_0$ be the one dimensional free module over $\cE$ with zero hermitian form.
Let $\tilde{K} = K \oplus \cE_0$.
Let us write the elements of $\tilde{K}$ in the form $(y,m)$ with $y \in K$ and
$m \in \cE$. Let $\tilde{G}$ be the subgroup of $\Aut(\tilde{K})$ generated
by reflections in $\tilde{\Phi} = \lbrace (x,m) \colon x \in \Phi, m \in \cE \rbrace$. 
Note that $\phi_{(x,m)}^{\omega} \in \tilde{G}$ if and only if $\phi_x^{\omega} \in G$. 
Consider the semi-direct product $K \rtimes G$,
in which the product is defined by $(x,g).(y,h) = (x + gy, gh)$. The faithful
action of $K \rtimes G$ on $K$ via affine transformations is given by
$(x, g) y = x + g y$.
\end{topic}
\begin{lemma}
Given the setup above, assume that $\abs{r}^2 = 3$ for all $r \in \Phi$ and $K' \supseteq p^{-1} K$.
\par
(a) The affine reflection group $\tilde{G}$ is isomorphic to the semi-direct product
$K \rtimes G $.
\par
(b) Let $r_0$ be a root of $\Phi$ such that the orbit $G r_0$ spans $K$ as a $\ZZ$--module.
If $\phi^{\omega}_{r_1}, \dotsb, \phi^{\omega}_{r_k}$
generate $G$, then $\phi^{\omega}_{(r_1,0)}, \dotsb, \phi^{\omega}_{(r_k,0)}$,
together with $\phi^{\omega}_{(r_0, 1)}$ generate $\tilde{G}$.
\par
(c) If $\phi^{\omega}_{r_0}$ commutes (resp. braids) with $\phi^{\omega}_{r_j}$, then
$\phi^{\omega}_{(r_0, 1)}$ commutes (resp. braids) with $\phi^{\omega}_{(r_j,0)}$. 
\label{lemma-affine}
\end{lemma}
\begin{proof}
(a) Identify $K$ (resp. $G$) inside $\tilde{K}$ (resp. $\tilde{G}$)
via $y \mapsto (y,0)$ (resp. $\phi^{\omega}_x \mapsto \phi^{\omega}_{(x,0)}$ ). 
Recall $p = 2 + \omega$. For $x \in K$, define 
\begin{equation}
t_x(y,n) = (y, n - p^{-1} \ip{x}{y} ).
\label{tx}
\end{equation}
The automorphisms $t_x$ of $\tilde{K}$ are called translations.
The subgroup of $\Aut(\tilde{K})$ generated by translations is isomorphic to
the additive group of $K$. For any root $(x,m) \in \tilde{\Phi}$, one has,
\begin{equation*}
 \phi_{(x,m)}^{\omega}(y,n) = (\phi_x^{\omega}(y), n - m p^{-1} \ip{x}{y} ) 
= \phi_{(x,0)}^{\omega} \circ t_{\bar{m} x} (y,n).
\end{equation*}
Let us write $\phi^{\omega}_{(x,0)} = \phi^{\omega}_x$. From the above equation, one has, in particular,
\begin{equation} 
 t_x =    (\phi^{\omega}_x)^{-1}  \circ \phi^{\omega}_{(x,1)} .
\label{eq-translation}
\end{equation}
From equation \eqref{eq-translation}, it follows that $ \phi^{\omega}_a t_x (\phi^{\omega}_a)^{-1} = t_{ \phi^{\omega}_a(x)}$.
So $(x, g) \mapsto t_x \circ g$ is an isomorphism from $ K \rtimes G$ onto $\tilde{G}$.
\par
(b) Let $G_1$ be the subgroup of $\tilde{G}$ generated by 
$\phi^{\omega}_{r_1}, \dotsb, \phi^{\omega}_{r_k}$ and $\phi^{\omega}_{(r_0,1)}$. Then $G \subseteq G_1$.
Let $x = g r_0$ be a root in the $G$--orbit of $r_0$.
Then $t_x = (\phi^{\omega}_x)^{-1} \phi^{\omega}_{(x,1)}  = (\phi^{\omega}_x)^{-1} g \phi^{\omega}_{(r_0,1)} g^{-1} \in G_1$.
Since the roots in the $G$--orbit of $r_0$ span $K$ as a $\ZZ$--module, it follows that $t_x \in G_1$
for all $x \in K$. The translations, together with $G$, generate $\tilde{G}$. So $\tilde{G} = G_1$.
\par
Part (c) follows from \ref{braid-commute} since $\lbrace (r_0,1), (r_j ,0)\rbrace$ and
$\lbrace r_0,r_j \rbrace$ have same gram matrix. 
\end{proof}
\begin{topic}{\bf Method to get an affine diagram:}
Lemma \ref{lemma-affine} applies to the root systems $\Phi(G_4)$, $\Phi(G_{25})$
and $\Phi(G_{32})$ and the corresponding lattices $D_4^{\cE}$, $E_6^{\cE}$ and $E_8^{\cE}$.
Similar result holds for the root systems $\Phi(G_{33})$ and $\Phi(G_{34})$ and the
corresponding lattices $K_{10}^{\cE}$ and  $K_{12}^{\cE}$, if one replaces
order three reflections by order two reflections and $p$ by $-\omega$.
In each of these cases, any root can be chosen as $r_0$ in part (b) of \ref{lemma-affine}.
\par
Further modifications are necessary for $G_{26}$.
In this case $p^{-1} K$ is not a subset of $K'$ but there is a sub-lattice $E_6^{\cE} \subseteq K$
such that $p^{-1} E_6^{\cE} \subseteq K'$. Accordingly, 
the translations $t_x$, given in \eqref{tx}, define automorphisms of $\tilde{K}$ only for
$x \in E_6^{\cE}$. The conclusion in part (a) is that $\tilde{G} \simeq E_6^{\cE} \rtimes G$. 
In part (b), any root of an order $3$ reflection can be chosen as $r_0$. The details are omitted.
\par
Let $G \in \lbrace G_4, G_{25}, G_{26}, G_{32}, G_{33}, G_{34} \rbrace$.
We take the diagram for $G$ obtained by method \ref{method} and
extend it by adding an extra node corresponding to a suitable root, thus
obtaining an affine diagram.
These are shown in figure \ref{fig-affine}. The extending node
is joined with dotted lines. The vertices of an affine diagram of type $G$
correspond to a minimal set of generators for the affine reflection group
$\tilde{G}$. The edges indicate the Coxeter relations among the generators.
Additional relations may be needed to obtain a presentation of $\tilde{G}$
(like those given in \ref{t-rel}).
\end{topic}
%
%
%
%
%
%
%
%
%
%
%
\appendix
\section{Root systems for some unitary reflection groups}
\label{ap-1}
We describe a root system for each unitary reflection group considered in section \ref{secweylvector}.
Notation: a set of co-ordinates marked with a line (resp. an arrow) above,
means that these co-ordinates can be permuted (resp. cyclically permuted).
\begin{topic}{$\mathbf{G_4, G_5, G_{25}, G_{26}, G_{32}}$:}
\label{a-1}
For each of these groups, a set of projective roots and the lattices spanned by
these roots are given in table \ref{table-1}.
The groups $G_5$, $G_{25}$ and $G_{32}$ are reflection groups of the $\cE$--root lattices
$D_4^{\cE}$, $E_6^{\cE}$ and $E_8^{\cE}$ respectively.
{\doublespacing
\begin{table}
\begin{tabular}{|c|c|c|c|}
\hline
$G$ & $\abs{\Phi_*(G)}$ & $\Phi_*(G)$ & $K$ \\ 
\hline
&&& \\
 $G_4$ &  $4$ & $(1,1,\omega^j)$, $(0,0,p)$   & $D_4^{\cE}$ \\
&&& \\
$G_5$ &  $4 + 4$ &   $\Phi_*(G_4)$ and $(p,p,0)$, $(1,1, -2 \omega^j)$  & $D_4^{\cE}$ \\
&&& \\
 $G_{25}$ & $12$ &$(\overline{p,0,0})$, $(1, \omega^j, \bar{\omega}^k)$ &   $E_6^{\cE}$ \\
&&& \\
$G_{26}$ &  $9+12$ &  $(\overline{1, -\omega^j, 0})$ and $\Phi_*(G_{25})$ &   \\
&&& \\
$G_{32}$ & $40$ & $(0,1,\omega^j, \bar{\omega}^k)$, $(1,\overrightarrow{\omega^j,-\bar{\omega}^k,0})$ & $E_8^{\cE}$  \\
&&& \\
\hline
\end{tabular}
\caption{}
\label{table-1}
\end{table}
}
\end{topic}
\begin{topic}{$\mathbf{G_{33}, G_{34}}$:}
\label{a-2}
These are reflection groups of the $\cE$--lattices $K_{10}^{\cE}$ and
$K_{12}^{\cE}$ respectively, where $K_{12}^{\cE}$ is a complex form of the Coxeter--Todd lattice
(see example 10a in chapter 7, section 8 of \cite{cs:splag}) and
$K_{10}^{\cE}$ is the orthogonal complement of any vector of minimal norm in $K_{12}^{\cE}$.
The minimal norm vectors of $K_{10}^{\cE}$ and $K_{12}^{\cE}$ form root
systems of type $G_{33}$ and $G_{34}$ respectively.
\end{topic}
\begin{topic}{$\mathbf{ G_8, G_{29}, G_{31} } $:}
\label{a-3}
These are the primitive unitary reflection groups defined over $\cG =\ZZ[i]$.
Let $q = 1 + i$. 
\par
The reflection group of the two dimensional $\cG$--lattice 
$D_4^{\cG} = \lbrace (x,y) \in \cG^2 : x + y \equiv 0 \bmod q \rbrace$ is $G_8$.
The six projective roots are $(q,0)$, $(0,q)$ and $(1, i^j)$. 
The reflection group contains order $4$ and order $2$ reflections in these
roots.
\par
A set of projective roots for $G_{29}$ can be chosen to be
\begin{equation*}
\Phi_*(G_{29}) = \lbrace (\overline{2,0,0,0}), (1, \overline{\pm 1,i,i} ), (1,\overline{\pm 1,-i,-i}), (\overline{q,\pm q,0,0}),
(1,\overline{\pm 1,i,-i}) \rbrace. 
\end{equation*}
There are a total of $4 + 6+ 6 + 12 + 12 = 40$ projective roots.
They span the $\cG$--lattice
$E_8^{\cG}  = \lbrace (x_1, \dotsb, x_4) \in \cG^4 \colon x_i \equiv x_j \bmod q, \sum_i x_i \equiv 0 \bmod 2 \rbrace $ whose real form is $E_8$. 
The minimal norm vectors of $E_8^{\cG}$ form a root system for $G_{31}$. 
The projective roots of $G_{31}$ can be chosen to be 
\begin{equation*}
\Phi_*(G_{31}) = \lbrace  (1, \pm 1 , \pm 1, \pm 1), (1, \overline{\mp 1, \pm 1, \pm 1}), (\overline{q, \pm i q,0,0}) \rbrace \cup \Phi_*(G_{29}).
\end{equation*}
There are $ 2 + 6 + 12 + 40 = 60$ projective roots.
\end{topic}
\begin{topic}{$\mathbf{G_{12}}:$}
\label{a-4}
following is a set of projective roots of $G_{12}$, defined over
$\ZZ[\sqrt{-2}]$:
\begin{equation*} 
\Phi_*(G_{12}) = \lbrace ( \overline{ 2, 0 }),  ( \overline{1,  \pm 1 \pm \sqrt{-2} }),  (\overline{ \sqrt{-2}, \pm \sqrt{-2}} ) \rbrace.
  \end{equation*}
  There are $2 + 8 + 2 = 12$ projective roots of norm $4$ and $G_{12}$ contains order $2$ reflections in these.
  \par
 The algorithm \ref{algo}, applied to $G_{12}$, yields three generators $s, t, v$ which satisfy the
relations $\lbrace  s t s v = t s v t , t s v s = v t s v \rbrace$ in each trial.
These relations are equivalent to the relations 
$ \lbrace s t u s = t u s t = u s t u \rbrace$ given in \cite{bmr:cr}
using the substitution $u  = c_s(v)$. So $\op{Braid}(G_{12})$ has a presentation given by
$\langle s, t, v \vert s t s v = t s v t , t s v s = v t s v  \rangle $.
\end{topic}
\begin{topic}{$\mathbf{G_{24}}:$}
\label{a-5}
The following is a set of projective roots of $G_{12}$, defined over $\ZZ[(1 + \sqrt{-7})/2]$:
\begin{equation*} 
\Phi_*(G_{24}) = \lbrace (\overline{2, 0, 0} ), ( \overline{ ( 1- \sqrt{-7})/2, \pm  (1 - \sqrt{-7})/2, 0} ) ,  (\overline{(1 + \sqrt{-7} )/2, \pm 1, \pm 1 })  \rbrace.
 \end{equation*}
The order $2$ reflections in these $(3+ 6 + 12) = 21$ roots of norm $4$ generate $G_{24}$. 
\par
The generators $s, t, u$ chosen by algorithm \ref{algo} satisfy the following relations: 
\begin{equation*} 
s t s = t s t, \; s t s t = u s u s , \; u t u t = t u t u, \; (s u t )^2 s = u (s u t)^2.
 \end{equation*}
 This happen to be the third presentation of  $\op{Braid}(G_{24})$
 given in \cite{bm:ep}. (note: the last relation above can be replaced by saying that $u$ braids with $s t s^{-1}$).
\end{topic}
\begin{topic}{$\mathbf{G_6, G_9}:$}
\label{a-6}
Let $\zeta_m = e^{2 \pi  i/m}$. The projective roots of $G_6$ (defined over $\ZZ[\zeta_{12}]$) can be chosen
to be 
\begin{equation*}
\Phi_*(G_{6})  = \lbrace (\overline{2,0}), (\pm \zeta_{12} q , q), ( \pm \omega q ,  q) \rbrace
\cup \lbrace (\pm 1, 1 + \zeta_{12}), ( 1 + \bar{\zeta}_{12} ,  \pm i) \rbrace.
\end{equation*} 
The group $G_6$ contains order $2$ reflections in the six projective roots of norm $4$
and order $3$ reflections in the four projective roots of norm $(3 + \sqrt{3})$.
\par
The projective roots of $G_9$ (defined over $\ZZ[\zeta_8]$) can
be chosen to be 
\begin{equation*}
\Phi_*(G_{9}) =  \Phi_*( G_{12}) \cup \lbrace (\pm 1, 1 + \sqrt{2}), (1 + \zeta_8, \pm i ( 1 + \zeta_8))\rbrace 
\end{equation*}
There are $12$ projective roots of norm $4$ (those of $G_{12}$) and six of
norm $(4 + 2 \sqrt{2})$.  The group $G_9$ contains order $2$ reflections
in all the roots and order $4$ reflections in the roots of norm
$(4 + 2 \sqrt{2})$.
\end{topic}
\begin{topic}{$\mathbf{G(de, e, n)}$:}
\label{a-7}
Let $\zeta_m = e^{2 \pi i/m}$ and let
$e_j$ be the $j$-th unit vector in $\CC^n$. The projective roots of $G(de, e, n)$ can be chosen to be
$\lbrace e_j , e_j - \zeta_{de}^t e_k :  1 \leq j < k \leq n , 1 \leq t \leq de \rbrace$.
For a detailed study of these groups, see \cite{bmr:cr}.
\end{topic}
%
%
%
%
%
\section{Proofs of some statements in section \ref{s-newd}}
\label{ap-2}
\begin{topic}{\bf Notations: }
We adopt the following notations.
If $x$ and $y$ are elements in a group, we write $ \bar{x} = x^{-1}$,  $c_x(y) : = x y x^{-1}$. We write
 $x \braid y$ (resp. $x \comm y$) as an abbreviation for ``$x$ braids with $y$" (resp. ``$x$ commutes with $y$").
 \par
In the proof of proposition \ref{p-mutate}, say, while working in the group $B_{31}$,
instead of writing $\varphi(w) = c_{\bar{a}_3}(a_4)$, we simply write $w =c_{\bar{a}_3}(a_4) $.
Similar abuse of notation is used consistently because it significantly simplifies writing.
\end{topic}
\begin{topic}{\bf Proof of Lemma \ref{l-cyclic}.}
\label{pf-l-cyclic}
\begin{proof}
Work in $\op{Cox}(\tilde{A}_n, \infty)$. Let $h_n = (x_0 x_1 \dotsb x_n)$ and $y_n  = h_n x_0$. 
\par
{\bf Step 1.} Observe that 
\begin{equation} 
y_n x_n = x_0 \dotsb x_{n-1} \; x_n x_0 x_n  = x_0 \dotsb x_{n-1} \; x_0 x_n x_0 = x_0 x_1 x_0 \; x_2 \dotsb x_n \; x_0
= x_1 y_n.
\label{e-yx}
 \end{equation}
 Since $\rho^2(y_n) = x_1^{-1} \rho(y_n) x_2$, it follows that
 \begin{equation} 
 \rho(y_n) = y_n
  \iff   \rho^2(y_n) = x_1^{-1} y_n   x_2 = y_n x_n^{-1} x_2.
  \label{e-rhoy}
  \end{equation}
Suppose $\rel{n+2}{x_0, \dotsb, x_n} $ is cyclic.
Then $\rho(y_n) = y_n$ implies $\rho^2(y_n) = y_n$, hence $x_n = x_2$ (from \eqref{e-rhoy}).
Conversely, if $n = 2$, then $\rho(y_n) = y_n$ implies $\rho^2(y_n) = y_n$ (by \eqref{e-rhoy}). So $\rel{n+2}{x_0, \dotsb, x_n}$
is cyclic.
\par
{\bf Step 2.} Using \eqref{e-yx}, we have,
\begin{equation} 
h_n^2 x_{n-1} = y_n  \; x_1 \dotsb x_{n-2} \;  x_n x_{n-1} x_n = y_n x_n (x_1 \dotsb x_n)  = x_1 h_n^2.
\label{e-hnx}
 \end{equation}
It follows that
\begin{equation} 
\rho(h_n^2) = h_n^2
  \iff   \rho^2(h_n^2) = x_1^{-1} h_n^2   x_1 = h_n^2 x_{n-1}^{-1} x_1.
  \label{e-rhoh2}
 \end{equation}
As in Step 1, (but using \eqref{e-rhoh2} in place of \eqref{e-rhoy}) we conclude that $\rel{2n + 2}{ x_0, \dotsb, x_n}$ is cyclic if and only if $x_{n-1} = x_1$.
This proves part (b).
\par
(c) Let $z_n = h_n^2 x_0$. From \eqref{e-hnx}, it follows that
\begin{equation*} 
\rho(z_n) = z_n 
\iff \rho^2(z_n) = x_1^{-1} z_n x_2 = x_1^{-1} h_n^2 x_0 x_2 = h_n^2 x_{n-1}^{-1} x_0 x_2 = z_n x_{n-1}^{-1} x_2
 \end{equation*}
 where the last equality is obtained by commuting $x_{n-1}$ and $x_0$, which holds since $n > 2$. 
 Part (c) follows.
\end{proof}
\end{topic}
\begin{topic}{\bf Proof of Lemma \ref{l-rel}.}
\label{pf-l-rel}
\begin{proof}
(a) $c_x(y) \comm z \iff z x y \bar{x} = x y \bar{x} z = \bar{y} x y z$ or equivalently
$y z x y = x y z x$. Further $z \comm c_x(y) = c_{\bar{y}}(x) \iff c_y(z) \comm x$. 
If $x \braid z$, then $c_x(y) \comm z \iff y \comm c_{\bar{x}}(z) = c_z(x)$.
\par
(b) Let $u = c_x(y)$. Then $c_u(z) = c_{ x y \bar{x}}(z)  = c_{x y z}(x) $ and 
$c_{\bar{z}}(u) = c_{\bar{z} x}(y) = c_{\bar{z} \bar{y}}(x)$. So
\begin{equation*}
u \braid z \iff c_{u}(z) = c_{\bar{z}}(u) \iff c_{x y z} (x) = c_{\bar{z} \bar{y}}(x) \iff \rel{6}{x, y, z}.
\end{equation*}
Since $x$ braids with $y$ and $z$, one has
\begin{equation*} 
x \braid c_y(z) \iff c_{\bar{y}}(x) \braid z \iff c_x(y) \braid z \iff y \braid c_{\bar{x}}(z) \iff y \braid c_z(x).
 \end{equation*}
(c) Let $u  = c_{x y}(z)$. Then $c_{u}(w) = c_{x y z \bar{y} \bar{x} }(w)  = c_{x y z \bar{y} w}(x) 
= c_{x y z w \bar{y}}(x) = c_{x y z w x}(y)$; and $c_{\bar{w}}(u) =
c_{\bar{w} x y}(z) = c_{ \bar{w} x \bar{z}}(y) = c_{\bar{w} \bar{z} x}(y) = c_{\bar{w} \bar{z} \bar{y}}(x)$.
So
\begin{equation*}
u \braid w \iff c_{u}(w) = c_{\bar{w}}(u) \iff c_{x y z w x} (y) = c_{\bar{w} \bar{z} \bar{y}}(x) \iff \rel{9}{x, y, z, w}.
\qedhere
\end{equation*}
\end{proof}
\end{topic}
\begin{remark}
Note that Lemma \ref{l-cyclic} implies that the relations of the form $\rel{m}{x_0, \dotsb, x_n}$ considered
in Lemma \ref{l-rel} are cyclic if $x, y, \dotsb$ satisfy the Coxeter relations of $\op{Cox}(\tilde{A}_n, \infty)$.
\end{remark}
%
%
%
%
%
%
\begin{topic}{\bf Proof of proposition \ref{p-mutate} for $G_{33}$ and $G_{34}$: }
 \label{pf-g33}
\begin{proof} (a) Work in the group
generated by $a_1, a_2, a_3, a_4, a_5, s, t, u, v, w$ subject to the relations: 
\begin{equation} 
s = a_1, \text{\; \;} t =a_2, \text{\; \;} u = a_3, \text{\; \;} v = c_{a_3}(a_4) ,\text{\;\;} w = c_{a_2 a_1} (a_5).
\label{e-cv33}
 \end{equation}
We find a sequence $R_{33,0} \subseteq R_{33,1} \subseteq \dotsb \subseteq R_{33,10} = R_{33}$
such that $R_{33, j}$ is a set of relations in the alphabet $a_1, a_2, \dotsb$
and $R_{33,j}$ is obtained from $R_{33,j-1}$ by adding a
single relation $w_j(a_1, a_2, \dotsb)$. We find another sequence
$R'_{33,0} \subseteq R'_{33,1} \subseteq \dotsb \subseteq R'_{33,10} = R'_{33}$
such that $R'_{33, j}$ is a set of relations in the alphabet $s,t, \dotsb$
and $R'_{33,j}$ is obtained from $R'_{33,j-1}$ by adding a
single relation $w'_j(s,t, \dotsb)$.
\par
For each successive $j \geq 1$,
we assume $R_{33, j-1}$ (or equivalently $R'_{33,j-1}$) and show that the relation
$w_j$ is equivalent to the relation $w'_j$, hence $R_{33, j}$ is equivalent
to $R'_{33, j}$. After $10$ such steps we find that $R_{33}$ is equivalent to $R'_{33}$.
\par
Let $R'_{33, 0}$ (resp. $R_{33,0}$) be the relations in $R'_{33}$ (resp. $R_{33}$) that only involve $s,t,u$ (resp.  $a_1, a_2, a_3$).
Clearly $R'_{33,0}$ and $R_{33,0}$ are equivalent under \eqref{e-cv33}.
Assume $R'_{33,0}$, or equivalently $R_{33, 0}$. The chain of equivalences below correspond to $j = 1, 2, 3, \dotsb$.
In the following chain of equivalences it is understood that the leftmost relation is $w_j$ and the rightmost one is $w'_j$.
The implicit assumptions in the $j$-th step are the relations $R_{33, j-1}$ and $R'_{33, j-1}$.  
\begin{equation*}
j = 1:  \;\;\;\; a_4 \comm a_1   \iff c_{\bar{u}}(v) \comm s \iff   v \comm s. 
\end{equation*}
\begin{equation*} 
j = 2: \;\;\;\; a_4 \braid a_3 \iff  c_{\bar{u}}(v) \braid u \iff v \braid u. 
 \end{equation*}
\begin{equation*}
j = 3: \;\;\;\; \rel{4}{a_3,a_4,a_2} \iff  c_{a_3}(a_4) \comm a_2 \iff v \comm t.
\end{equation*}
\begin{equation*}
j = 4: \;\;\;\; a_4 \braid a_2 \iff c_{\bar{u}}(v) \braid t \iff v \braid c_u(t) = c_{\bar{t}}(u) \iff v \braid u.
\end{equation*}
Note that $R'_{33,4}$ is the set of relations in $R'_{33}$ 
not involving $w$ and $R_{33,4}$ is the set of relations in $R_{33}$ not involving $a_5$.
At this stage we have shown that $R'_{33,4}$ is equivalent to $R_{33,4}$. 
 These are our implicit assumptions in the next step given below:
\begin{equation*} 
a_5 \comm a_2 \iff c_{\bar{s} \bar{t} }(w) \comm t \iff w \comm c_{t s }(t)  \iff w \comm s. 
 \end{equation*}
\newline
\begin{equation*} 
a_5 \braid a_1 \iff c_{\bar{s} \bar{t}} (w) \braid s \iff w \braid c_t(s) = c_{\bar{s}}(t) \iff w \braid t.
 \end{equation*}
\newline
\begin{equation*} 
a_5 \braid a_3 \iff c_{\bar{s} \bar{t}}(w) \braid u \iff w \braid c_{t}(u) \iff \rel{6}{t, u, w}.
 \end{equation*}
 \newline
 \begin{equation*} 
 \rel{9}{a_3, a_2, a_1, a_5} \iff a_5 \braid  c_{a_3 a_2}(a_1) 
 \iff c_{\bar{s} \bar{t}}(w) \braid c_{u t } (s) = c_{ u \bar{s}}(t)
= c_{\bar{s} \bar{t}}(u) \iff w \braid u.
 \end{equation*}
\newline
\begin{equation*} 
\rel{4}{a_3, a_4, a_5} \iff a_5 \comm c_{a_3}(a_4) \iff c_{\bar{s} \bar{t}}(w) \comm v \iff w \comm v.
 \end{equation*}
\newline
\begin{align*} 
a_5 \braid a_4 
\iff c_{\bar{s} \bar{t}} (w) \braid c_{\bar{u}}(v) 
 \iff c_{\bar{t}}(w) \braid c_v(u)
 \iff c_{\bar{t}}(w)  \braid u
 &\iff w \braid  c_{t}(u) \\ 
&\iff \rel{6}{t, u, w}.  
 \end{align*}
This proves \ref{p-mutate} for $G_{33}$.
\par
(b) {\bf Step 1: } Assume the relations $R_{34}$. Let $s, t, u, v, w$ be as in \eqref{e-cv33}. Let
 \begin{equation} 
 x = c_{a_4 a_2 a_1 a_5 a_3}(a_6).
 \label{e-cv34}
  \end{equation}
  We have to check that $\lbrace s, t, u,v,w,x \rbrace$ satisfy $R'_{34}$. Part (a) implies that $\lbrace s, t, u,v,w\rbrace$
  satisfy the relations in $R'_{34}$ that does not involve $x$. The following observation is helpful in checking the
relations involving $x$.
Let \begin{equation*} 
\alpha = c_{a_{5} a_3}(a_6), \text{\; \;}  \delta =\bar{a}_1 \bar{a}_2 \bar{a}_4  , \text{\; so that \;} x = c_{\bar{\delta}}(\alpha).
 \end{equation*}
{\it Claim: $\alpha$ commutes with $a_2$, $a_3$ and $a_4$ and braids with $a_5$.}
\begin{proof}[proof of claim]
\begin{equation*}
a_2 \comm c_{a_5 a_3}(a_6) \iff a_2 \comm c_{a_3}(a_6) \iff \rel{4}{a_2, a_3, a_6}. 
\end{equation*}
\begin{equation*} 
c_{a_5 a_3}(a_6) \comm  a_3 \iff c_{a_3}(a_6) \comm c_{\bar{a_5}}(a_3) = c_{a_3}(a_5) \iff a_6 \comm a_5. 
 \end{equation*}
\begin{equation*} 
a_4 \comm c_{a_5 a_3}(a_6) = c_{a_5 \bar{a}_6}(a_3) = c_{\bar{a}_6 a_5}(a_3) \iff a_4 \comm c_{a_5}(a_3) \iff \rel{4}{a_3, a_4, a_5}. 
\end{equation*}
\begin{equation*}
a_5 \braid c_{ a_5 a_3 }(a_6) \iff a_5 \braid c_{a_3}(a_6) = c_{\bar{a}_6}(a_3) \iff a_5 \braid a_3. \qedhere
 \end{equation*}
\end{proof}
Now we can check the relations in $R_{34}'$ involving $x$. Note that
\begin{equation*} 
c_{\delta}(a_1) = c_{ \bar{a}_1 \bar{a}_2}(a_1) = a_2,  \; \; c_{\delta}(a_2) = c_{ \bar{a}_1}(a_4) = a_4, 
\; \; c_{\delta}(a_3) = c_{\bar{a}_1 \bar{a}_2  a_3}(a_4) = c_{\bar{a}_1 a_3}(a_4) =   c_{a_3}(a_4).
 \end{equation*}
Conjugating by $\delta$, we find that
$x$ commutes with $s = a_1$ (resp. $ t= a_2$ or $ u = a_3$) if and only if $\alpha$ commutes with 
$c_{\delta}(a_1) = a_2$ (resp. $c_{\delta}(a_2) = a_4$, or $c_{\delta}(a_3) = c_{a_3}(a_4)$).
\par	
Since $x$ and $a_3$ commutes, $x$ braids with $v = c_{a_3}(a_4)$ if and only $x$ braids with $a_4$. Note that
$c_{\delta}(a_4) = c_{\bar{a}_1 \bar{a}_2}(a_4) = c_{\bar{a}_1 a_4}(a_2) = c_{a_4 \bar{a}_1}(a_2) = c_{a_4 a_2}(a_1)$. It follows that
\begin{align*} 
x \braid v \iff x \braid a_4 \iff \alpha \braid c_{\delta}(a_4) =  c_{a_4 a_2}(a_1)   
& \iff \alpha = c_{a_5 a_3}(a_6) \braid a_1 \\
& \iff \rel{9}{a_5,a_3,a_6,a_1}.
 \end{align*} 
Since $x$ commutes with $a_1$ and $a_2$, one has $x \comm w = c_{a_2 a_1}(a_5)$ if and only if $x \comm a_5$. Using the
relations satisfied by $\alpha$, one obtains
$x= c_{ a_4 a_2 a_1}(\alpha) = c_{a_4 a_2 \bar{\alpha}}(a_1) = c_{\bar{\alpha} a_4 a_2}(a_1)$, so $c_{\alpha}(x) = c_{a_4 a_2}(a_1)$.
Also $c_{\alpha}(a_5) = c_{\bar{a}_5}(\alpha) = c_{a_3}(a_6) = c_{\bar{a}_6}(a_3) $. So
\begin{align*} 
a_5 \comm x  \iff c_{\bar{a}_6}(a_3) \comm c_{ a_4 a_2}(a_1) 
& \iff c_{\bar{a}_4}(a_3) \comm c_{a_6} (c_{a_2}(a_1))  = c_{a_2}(a_1)  \\
& \iff c_{\bar{a}_2}(c_{a_3}(a_4)) = c_{a_3}(a_4) \comm a_1. 
 \end{align*}
\par
{\bf Step 2: } Conversely, assume the relations $R'_{34}$. Let 
\begin{equation*} 
a_1 = s, \text{\; \;} a_2 = s, \text{\; \;} a_3 = u, \text{\; \;} a_4 = c_{\bar{u}}(v) ,\text{\;\;} a_5 = c_{\bar{s}\bar{t} }(w) ,
\text{\; \;} a_6 = c_{\bar{u} \bar{s} \bar{t} \bar{w} \bar{u} \bar{v}}(x).
 \end{equation*}
Let
\begin{equation*} 
\beta = c_{ v u}(w) , \; \;  \gamma = t s u \bar{x}.
 \end{equation*}
It is useful to note the following relations: 
\begin{equation} 
\text{$\beta$ commutes with $s$ and $u$, and $\beta$ braids with $c_{t}(s)$ and $c_{t}(u)$}.
\label{e-tuwcon}
 \end{equation}
\begin{proof}[proof of \eqref{e-tuwcon}] Since $s$ commutes with $u, v$ and $w$, it commutes with $\beta$.
Next, $ \beta = c_{v u}(w)\comm u$ if and only if  $w \comm c_{\bar{u} \bar{v} }(u) = v$.
Now, $\beta$ braids with $c_t(s) = c_{\bar{s}}(t)$ (resp. $c_t(u) = c_{\bar{u}}(t)$) if and only if $\beta$ braids with $t$.
Finally $\beta = c_{v u}(w) \braid t \iff c_{u}(w) \braid t \iff \rel{6}{t, u, w}$.
\end{proof}
Note that
\begin{equation*} 
\beta  = c_{ v \bar{w}}(u) = c_{\bar{w} v}(u) = c_{\bar{w} \bar{u}}(v) ; \text{ so \;} 
x = c_{ x \bar{u} \bar{s} \bar{t}  \; \bar{w} \bar{u} }(v)  = c_{\bar{\gamma} }(\beta).
 \end{equation*} 
Now we verify that $a_1, \dotsb, a_6$ satisfies the relations in $R_{34}$ involving $a_6$
 by reducing them to the relations in $R'_{34}$ or those mentioned in \eqref{e-tuwcon}.
\begin{equation*} 
a_6 \braid a_1 \overset{c_{\gamma} }{\iff} \beta \braid c_{t s u \bar{x}} ( s) = c_t(s).
 \end{equation*}
\newline
 \begin{equation*} 
 a_6 \braid a_2  \overset{c_{s\gamma }}{\iff}  c_{s }(\beta) \braid c_{s t s u \bar{x} }(t)   
\iff   \beta \braid  c_{s t s u}(t) = c_{s t s \bar{t}}(u) = c_{t s }(u) = c_t(u).
\end{equation*}
\newline
\begin{equation*} 
a_6 \braid a_3  \overset{c_{\gamma}}{\iff}  \beta \braid    c_{ t s u \bar{x}} (u) = c_{t}(u) .
\end{equation*}
\newline
\begin{equation*}
a_6 \comm a_4
\iff x \comm c_{ v u w s t u \; \bar{u}}(v) = c_{ v u w s t}(v) = c_{v u} (v) = u. 
\end{equation*} 
\newline
\begin{equation*} 
a_6 \comm a_5 
\overset{c_{\gamma}}{\iff}  c_{\bar{w} \bar{u} }(v)  \comm c_{t s u \bar{x}  \bar{s} \bar{t}}(w)  
= c_{ t u \bar{t}}(w) = c_{ c_t(u) }(w) = c_{\bar{w}}(c_t(u)) = c_{\bar{w}\bar{u}}(t) 
\iff  v \comm t.
 \end{equation*}
 \newline
  \begin{equation*} 
 \rel{4}{a_2, a_3, a_6}
  \iff  a_6 \comm c_{a_2}(a_3) = c_t(u) \overset{c_{\gamma}}{\iff}  \beta \comm c_{t s u \bar{x} t }( u)   =   c_{t s }(t) = s.
  \end{equation*}
\newline
 \begin{equation*} 
 \rel{4}{a_2, a_1, a_6} 
 \iff  a_6 \comm  c_{a_2}(a_1) = c_{t}(s)  
\overset{c_{\gamma}}{\iff}   \beta \comm c_{t u s t }(s) = c_{t u}(t) = u.
\end{equation*}
  \newline
 \begin{align*} 
 \rel{9}{a_5, a_3, a_6, a_1} 
 \iff c_{ a_5 a_3}(a_6) \braid a_1 
 \iff & a_6 \braid c_{\bar{a}_3 \bar{a}_5}(a_1)  = c_{ \bar{u} \bar{s} \bar{t} \bar{w} t s}(s)  = c_{\bar{u} \bar{s} \bar{t} \bar{w} \bar{s}}(t) \\
 \iff & c_{\bar{u} \bar{v}} (x) \braid c_{\bar{s}}(t) \iff c_{x v}(u) \braid t \iff u \braid t.
  \end{align*}
 \end{proof}
 \end{topic}
%
%
%
%
%
%
\begin{topic}{\bf Proof of proposition \ref{p-mutate} for $G_{29}$.}
\label{pf-g29}
\begin{proof}
{\bf Step 1:}
Define $s, \dot{s}, t, v, u$ in $B_{29}$ by
\begin{equation*} 
s = c_{a_4 a_2 a_3}(a_1), \;\; \dot{s} = c_{ \bar{a}_1 \bar{a}_3 \bar{a}_2 } (a_4), \;\;  t = a_1 , \;\; v = a_2,  \;\; u = a_3.
\end{equation*} 
One has
 \begin{align*} 
 c_{\bar{a}_4 \bar{a}_1 \bar{a}_3 \bar{a}_2} (a_4) 
 & = c_{ a_3 \bar{a}_1  \bar{a}_3 \bar{a}_4 \bar{a}_1 \bar{a}_2} (a_4) \text{\;\; 
 (using $\rel{4}{a_3, a_1, a_4}$ in the form $\bar{a}_4 \bar{a}_1 \bar{a}_3  = a_3 \bar{a}_1  \bar{a}_3 \bar{a}_4 \bar{a}_1 $    )    }  \\
 & = c_{a_3 \bar{a}_1 \bar{a}_3} (a_2)  \text{\;\; (using $\rel{4}{a_4, a_1, a_2}$  in the form $c_{\bar{a}_4 \bar{a}_1 \bar{a}_2}(a_4) = a_2 $ )} \\
 &= c_{a_2 a_3} (a_1)  \text{\;\; (using $\rel{6}{a_3,a_1, a_2}$ in the form $a_2 \braid c_{a_3}(a_1)$   )}, 
 \end{align*}
that is,
\begin{equation*} 
s = \dot{s}.
 \end{equation*}
 Now we verify that $s, t, v, u$ satisfy  $R'_{29}$. Only the relations involving $s$ require some work.
 \begin{align*} 
 s \comm v \iff c_{ a_4 a_2 a_3 }(a_1) \comm a_2 \iff c_{a_2 a_3}(a_1) \comm c_{\bar{a}_4}(a_2) = c_{a_2}(a_4) 
 & \iff  c_{a_3}(a_1) \comm a_4 \\
 & \iff \rel{4}{a_3, a_1, a_4}.
 \end{align*}
 \begin{align*} 
 s \comm u \iff c_{ \bar{a}_1 \bar{a}_3 \bar{a}_2 }(a_4) \comm a_3 \iff  c_{ \bar{a}_3 \bar{a}_2 }(a_4) \comm c_{a_1}(a_3) = c_{\bar{a}_3}(a_1) 
 & \iff c_{ \bar{a}_2 }(a_4) \comm  a_1 \\
 & \iff \rel{4}{ a_4, a_2, a_1}.
 \end{align*}
\begin{equation*} 
c_{\dot{s}}(a_1) = c_{\bar{a}_1 \bar{a}_3 \bar{a}_2 a_4 a_2 a_3 a_1}(a_1) 
= c_{\bar{a}_1 \bar{a}_3 \bar{a}_2 a_4 a_2 a_3 }(1)  = c_{\bar{a}_1 \bar{a}_3 \bar{a}_2}(s)
=c_{\bar{a}_1}(s); 
 \end{equation*}
where the last equality uses $s \comm a_2$ and $s \comm a_3$. So $s \braid t$.
\par
{\bf Step 2:} Conversely, define $a_1, a_2, a_3, a_4 \in B'_{29}$, by 
\begin{equation*} 
a_1 = t, \;\; a_2 = v, \;\; a_3 = u, \text{\; and \;} a_4 = c_{ v u t} (s) = c_{v u \bar{s}}(t) = c_{\bar{s} v u}(t).
 \end{equation*}
 We have to check that $a_1, a_2, a_3, a_4$ satisfies $R_{29}$. 
 Only the relations involving $a_4$ require some verification, which is done below.
(We write $t \braid_4 v$ to denote the relation $\rel{4}{t, v}$).
\begin{equation*} 
a_4 \braid a_2 \iff c_{\bar{s} v u}(t) \braid v \iff c_u(t) \braid v \iff \rel{6}{u, t, v}.
 \end{equation*}
 \begin{equation*} 
 a_4 \braid_4 a_3 \iff c_{\bar{s} v u}(t) \braid_4 u \iff t \braid_4 c_{\bar{u} \bar{v} }(u)  = v.
  \end{equation*}
 \begin{align*} 
 a_4 \braid_4 a_1 & \iff c_s(t) \braid_4 c_{v u}(t) = c_{c_u(\bar{t}) }(v)  \text{\; \; (since $c_u(t) \braid v$)} \\  
 & \iff c_{\bar{t}}(s) \braid_4 c_{ \bar{t} \bar{u} t}(v) \iff 
s \braid_4 c_t(v) \iff c_s(t) = c_{\bar{t}}(s) \braid_4 v \iff t \braid_4 v.
\end{align*}
\end{proof}
\end{topic}
\begin{topic}{\bf Proof of proposition \ref{p-mutate}  for $G_{31}$: }
\label{pf-g31}
\begin{proof}
{\bf Step 1:} Assume $a_1, \dotsb, a_5$ satisfy $R_{31}$. 
 Define $s, t, u,v,w$ by
 \begin{equation} 
s = a_5, \;\; t = a_1, \;\; u = a_3, \;\; w = c_{\bar{a}_3}(a_4),  \;\; v = c_{a_4}(a_2).
\label{e-cv31}
 \end{equation}
 One has to check that $s, t, \dotsb$ satisfy $R'_{31}$.
The braid relations between $s, t, u$ in $R'_{31}$ are the same as the braid relations between $a_5, a_1, a_3$ in $R_{31}$.
The rest of $R'_{31}$ is verified below.
\begin{equation*} 
t \comm w \iff 1 \comm c_{\bar{a}_3}(a_4) \iff c_{a_3}(a_1) \comm a_4 \iff \rel{4}{a_3,a_1,a_4}.
 \end{equation*}
\begin{equation*} 
s \comm v \iff a_5 \comm c_{a_4}(a_2) \iff \rel{4}{a_4,a_2,a_5}.
 \end{equation*}
\begin{equation*} 
t \comm v \iff a_1 \comm c_{a_4}(a_2) \iff \rel{4}{a_4,a_2,a_1}.
 \end{equation*}
 \begin{equation*} 
 u \braid v \iff a_3 \braid c_{a_4}(a_2) \iff \rel{6}{a_2,a_3,a_4} \text{\; \; (see lemma \ref{l-rel}(b))}.
  \end{equation*}
\begin{align*} 
v \braid w \iff c_{a_4}(a_2) \braid c_{\bar{a}_3}(a_4) \iff c_{a_3 a_4}(a_2) \braid a_4 & 
\iff c_{c_{a_4}(\bar{a}_2)} (a_3) \braid a_4 \text{\;(using $\rel{6}{a_2,a_3,a_4}$)}\\
 &\iff a_3 \braid c_{ a_4 a_2 \bar{a}_4}(a_4) = c_{a_4 a_2}(a_4) = a_2.
\end{align*}
Finally $s u w = u w s = w s u$ follows from $\rel{3}{a_5, a_4, a_3}$ and $\rel{3}{a_4, a_3, a_5}$ as shown below:
\begin{align*} 
& s u w = a_5 a_3 c_{\bar{a}_3}(a_4) = a_5 a_4 a_3, \;\; u w s = a_3 c_{\bar{a}_3}(a_4) a_5 
= a_4 a_3 a_5, \\ 
& w s u =  \bar{a}_3 \; a_4 a_3 a_5 \; a_3 = \bar{a}_3 a_3 a_5 a_4 a_3 = a_5 a_4 a_3.
 \end{align*}
 {\bf Step 2: } Conversely assume that $s,t,u,v,w$ satisfies $R'_{31}$. Invert the equations in \eqref{e-cv31} to obtain
 \begin{equation*}
 a_5 = s,\;\;  a_1 = t, \;\;  a_3 = u, \;\;  a_4 = c_u(w) = c_{\bar{s}}(w), \;\; a_2 = c_{ u \bar{w} \bar{u} }(v) = c_{\bar{s} \bar{w}}(v).
\end{equation*}
Again, The braid relations between $a_5, a_1, a_3$ are the same as the braid relations between $s, t, u$.
The rest of  $R_{31}$ is verified below.
\begin{equation*} 
a_3 \braid a_2 \iff u \braid c_{\bar{s} \bar{w}}(v) \iff c_{w s}(u) \braid v \iff u \braid v \text{\;\; (since $w s u = u w s$ implies $c_{w s}(u) = u$)}.
 \end{equation*}
\begin{equation*} 
 a_4 \braid a_2 \iff  c_{\bar{s}}(w) \braid c_{\bar{s} \bar{w}}(v) \iff w \braid v.
\end{equation*}
Now, since we know the braid relations between $a_1,\dotsb, a_5$, lemma \ref{l-rel} can be applied to conclude 
$\rel{4}{a_3,a_1,a_4} \iff c_{a_3}(a_1) \braid a_4$ etc. So the equivalences given in step 1, show that
$\rel{4}{a_3, a_1,a_4}, \rel{4}{a_4,a_2,a_5}, \rel{4}{a_4,a_2,a_1}$ and $\rel{6}{a_2,a_3,a_4}$ holds. One has
\begin{equation*} 
\rel{4}{a_1,a_5,a_2} \iff a_2 \comm c_{a_1}(a_5)  \iff  c_{\bar{s} \bar{w}}(v) \comm c_t(s) = c_{\bar{s}}(t)  \iff c_{\bar{w}}(v) \comm t.
 \end{equation*}
Finally
\begin{equation*} 
a_3 a_5 a_4 =  u s c_u(w) = u \; s u w \; \bar{u} = u w s u \bar{u} = u w s, \; \;
a_5 a_4 a_3 = s c_u(w) u = s u w, \;\;  a_4 a_3 a_5 = u w s.
 \end{equation*}
 \end{proof}
 \end{topic}
\begin{remark}
\label{r-remainingrel}
(1)  {\it  The relation $\rel{9}{ a_4, a_2, a_1, a_5}$ holds in $B_{33}$. }
\begin{proof} In the notation of \ref{pf-g33}, 
one has $\rel{9}{ a_4, a_2, a_1, a_5} \iff a_5 \braid c_{a_4 a_2}(a_1) \iff c_{\bar{s} \bar{t}}(w) \braid c_{\bar{u} v u t }(s) = c_{ v u \bar{v}  t}(s)
=c_{v u t}(s) = c_{v \bar{s} u}(t) = c_{v \bar{s} \bar{t}}(u) \iff w \braid u$.
\end{proof}
\par
(2) In step 2 of \ref{pf-g29}, we showed $a_1 \braid s$. So 
one has 
\begin{align*} 
a_1 \braid  c_{ \bar{a}_1 \bar{a}_3 \bar{a}_2}(a_4) \iff a_1 \braid c_{\bar{a}_3 \bar{a}_2}(a_4) 
& \iff c_{\bar{a}_1}(a_3) = c_{a_3}(a_1) \braid  c_{a_4}(a_2) \\
& \iff a_3 \braid c_{a_4}(a_2) \text{\;(using $a_1 \comm c_{a_4}(a_2)$)}.
 \end{align*}
So {\it $\rel{6}{a_2, a_3, a_4}$ holds in $B_{29}$}. It follows that 
{\it there is an automorphism $\sigma$ of $B_{29}$ such that $\sigma^2 = 1$, 
$\sigma(a_1) = a_4^{-1}$, $\sigma(a_2) = a_3^{-1}$ (a lift of the diagram automorphism
of $D_{29}$).}
\par
Let $B_{29,+}$ be the monoid generated by $a_1, a_2, a_3, a_4$ subject to the relations $R_{29}$. 
Then $\tilde{\sigma}(x) = \sigma(x)^{-1}$ is an anti-automorphism of the monoid $B_{29,+}$ of order $2$.
To apply
 $\tilde{\sigma}$ to a word representing an element of $B_{29,+}$ one interchanges the occurrences of $a_1$ and
 $a_4$ (resp. $a_2$ and $a_3$) and writes the word from right to left. In other words, the diagram automorphism
 of $D_{29}$ lifts to an automorphism from $B_{29,+}$ to its opposite monoid $B_{29,+}^{op}$.
 \par
 Observe that $\dot{s} = s$ translates into $\sigma(s) = s$.
 The verification of $s \comm u$ in \ref{pf-g29} is just application of $\sigma$ to the verification of $s \comm v$.
 \end{remark}
%
%
%
%
%
%
%
%

\end{document}